\documentclass[11pt]{article}
\usepackage[english]{babel}
\usepackage{indentfirst}
\usepackage{latexsym}
\usepackage[usenames,dvipsnames]{color}
\usepackage[colorlinks=true, bookmarksnumbered=true, bookmarksopen=true,
bookmarksopenlevel=3, pdfstartview=FitH, linkcolor=cyan, pdfmenubar=true,
pdftoolbar=true, bookmarks=true,citecolor=cyan, urlcolor=magenta,
filecolor=cyan,menucolor=black,plainpages=false,pdfpagelabels]{hyperref}
\usepackage[lmargin=2cm,tmargin=2cm,bmargin=2cm,rmargin=2cm]{geometry}
\usepackage{amsbsy,amsmath,amsfonts,amssymb,amsthm}
\usepackage{times}
\usepackage{graphicx}
\usepackage[utf8]{inputenc}
\usepackage[T1]{fontenc}
\usepackage{amsmath}
\usepackage{amsfonts}
\usepackage{amssymb}

\numberwithin{equation}{section}
\newtheorem{prop}{Proposition}[section]
\newtheorem{definition}[prop]{Definition}
\newtheorem{theorem}[prop]{Theorem}

\newtheorem{remark}[prop]{Remark}
\newtheorem{lemma}[prop]{Lemma}

\newtheorem{proposition}[prop]{Proposition}

\def\fin { \vskip 0pt \hfill $\diamond$ \vskip 12pt}

\begin{document}

\title{On the well-posedness in Besov-Herz spaces for the inhomogeneous \\
incompressible Euler equations}
\author{{Lucas C. F. Ferreira$^{1}$}{\thanks{LCFF was partially supported by FAPESP (Grant: 2020/05618-6) and CNPq (Grant: 308799/2019-4), BR. Email: lcff@ime.unicamp.br (Corresponding author).}}, \ Daniel F. Machado$^{2}${\thanks {DFM was supported by CAPES (Finance Code 001), BR. Email: daniellmath@gmail.com.}}\\{\small $^{1,2}$ Universidade Estadual de Campinas, IMECC-Departamento de
Matem\'{a}tica,}\\{\small Rua S\'{e}rgio Buarque de Holanda, CEP 13083-859, Campinas, SP,
Brazil.}}
\date{}
\maketitle

\begin{abstract}
In this paper we study the inhomogeneous incompressible Euler equations in
the whole space $\mathbb{R}^n$ with $n\geq3$. We obtain well-posedness and
blow-up results in a new framework for inhomogeneous fluids, more precisely
Besov-Herz spaces that are Besov spaces based on Herz ones, covering
particularly critical cases of the regularity. Comparing with previous works
on Besov spaces, our results provide a larger initial data class for a
well-defined flow. For that, we need to obtain suitable linear estimates for
some conservation-law models in our setting such as transport equations and
the linearized inhomogeneous Euler system.

\medskip


\end{abstract}


\renewcommand{\abstractname}{Abstract}


\section{Introduction}

We consider the density-dependent incompressible Euler equations
\begin{equation} \label{sist:Euler_0}
\left\{
\begin{array}{l}
\partial _{t}\rho +u\cdot \nabla \rho =0, \\
\rho \left( \partial _{t}u+u\cdot \nabla u\right) +\nabla \pi =\rho f, \\
\mathrm{div}\;u=0, \\
(\rho ,u)_{|_{t=0}}=(\rho _{0},u_{0}),
\end{array}
\right. \hspace{0.5cm} (x,t)\in \mathbb{R}^{n}\times \mathbb{R}^{+},
\end{equation}
where $n\geq 3,$ $\rho $ is the density, $u=(u_{1},\dots ,u_{n})$ is the velocity field of the fluid and $\pi $ is the scalar pressure. Moreover, $f$ denotes a time-dependent vector field representing a given external force and we assume that there are three constants $\underline{\rho },\overline{\rho },\widetilde{\rho }>0$ such that
\begin{equation} \label{cond-basic-density-1}
0<\underline{\rho }\leq \rho _{0}(x)\leq \overline{\rho }<\infty
\end{equation}
and $\rho _{0}(x)\rightarrow \widetilde{\rho }$ when $|x|\rightarrow \infty $. For the sake of simplicity, spaces of scalar and vector functions are denoted in the same way. More precisely, for a Banach space $X$ and $u=(u_{1},\dots ,u_{n})$, the vector $u\in X$ means that $u_{1},\dots,u_{n}\in X$.

In comparison with the standard incompressible Euler equations (i.e., with constant density), see e.g. \cite{Bourgain_2015, Chae_2003_3, Chae-2004, Chemin_2, Lucas_Herz, Takada, Vishik-1998} and their references, there are relatively few works concerning to the well-posedness (existence and uniqueness) for (\ref{sist:Euler_0}) in the nonconstant density case. In fact, the well-posedness analysis for (\ref{sist:Euler_0}) in Sobolev and Besov type spaces is involved, especially in lower regularity spaces. In bounded domains of $\mathbb{R}^{2}$ or $\mathbb{R}^{3}$, local well-posedness results were obtained by Beir\~{a}o da Veiga and Valli \cite{Beirao-Valli-80-1, Beirao-Valli-80-2, Beirao-Valli-80-3} by considering initial data of class $C^{\infty }$ or with high H\"{o}lder-type regularity. Later, still in bounded domains $\Omega \subset \mathbb{R}^{n}$ but also treating high dimensions, Valli and Zaj\k{a}czkowski \cite{Valli-Za-88} proved the well-posedness for (\ref{sist:Euler_0}) in the Sobolev space $W^{2,p}(\Omega ).$ This setting was also employed by Itoh and Tani \cite{Itoh-99} with $\Omega $ being either a bounded or unbounded smooth domain $\Omega $ of $\mathbb{R}^{3}$. In the works \cite{Itoh-94} and \cite{Itoh-95}, Itoh considered the whole space $\Omega =\mathbb{R}^{3}$ and obtained local existence of solution with initial data $(\rho _{0}-c,u_{0})\in $ $H^{3}\times H^{3}$ or $(\nabla \rho _{0},u_{0})\in $ $H^{2}\times H^{3}$, respectively, where $c>0$ is a constant. In \cite{Danchin-2006}, the results of \cite{Itoh-95} were extended to $H^{s}(\mathbb{R}^{n})$\ with $s>1+n/2$ and $n\geq 2.$ Chae and Lee \cite{Chae_2003_2} treated the critical regularity case by analyzing the local well-posedness in the nonhomogeneous Besov space $B_{p,q}^{s}$ with $s=1+n/2,$ $p=2$ and $q=1$ (see also \cite[Section 7]{Danchin-2006}). In \cite{Danchin_2010}, Danchin proved the local well-posedness by considering $u_{0}\in L^{2}\cap B_{p,q}^{s}$ and $\rho_{0} $ is as in (\ref{cond-basic-density-1}) such that $\nabla \rho _{0}\in B_{p,q}^{s-1}$, where $1<p<\infty $ and $s>1+n/p$ with $1\leq q\leq \infty $ or $s=1+n/p$ with $q=1.$ The endpoint cases $p=1$ and $p=\infty $ were treated posteriorly by Danchin and Fanelli \cite{Danchin-2011}. The infinite energy case, even in a localized sense, are more subtle and the results found in the literature of (\ref{sist:Euler_0}) require additionally a smallness condition on the initial density (say, small perturbation of a constant) due to the influence of the pressure term. In this direction, we have the works \cite{Zhou_2010} and \cite{Wei_2013} where the local well-posedness for (\ref{sist:Euler_0}) was obtained in $B_{p,q}^{s}$ by assuming $s>1+n/p$ with $1\leq q\leq \infty $ and $s=1+n/p$ with $q=1,$ respectively. Moreover, let us point out that all aforementioned results considered (explicitly or implicitly) $\rho _{0}$ being bounded and bounded away from zero, that is, the basic condition (\ref{cond-basic-density-1}) which is preserved by the evolution of the density $\rho $.

Making the change $a=1/\rho -1$, we can rewrite (\ref{sist:Euler_0}) in the equivalent form
\begin{equation} \label{sist:Euler_1}
\left\{ \begin{array}{l}
\partial _{t}a+u\cdot \nabla a=0, \\
\partial _{t}u+u\cdot \nabla u+(1+a)\nabla \pi =f, \\
\mathrm{div}\;u=0, \\
(a,u)_{|_{t=0}}=(a_{0},u_{0}),
\end{array}
\right. \hspace{0.5cm} (x,t)\in \mathbb{R}^{n}\times \mathbb{R}^{+},
\end{equation}
where we still employ the word \textquotedblleft density\textquotedblright\ to refer to the unknown $a$.

In this work we analyze the local well-posedness of (\ref{sist:Euler_1}) in a new framework of Besov type for inhomogeneous fluids, namely nonhomogeneous Besov-Herz spaces $BK_{p,q,r}^{\alpha ,s}$ which are nonhomogeneous Besov spaces based on Herz spaces $K_{p,q}^{\alpha }$. Our results cover the critical case of the regularity as well as values above it. Also, we provide a blow-up criterion. The Herz $K_{p,q}^{\alpha }$ spaces were introduced initially by Herz \cite{Herz} as an appropriate environment for Fourier transform action. After, a new characterization for the $K_{p,q}^{\alpha }$-norm based on the $L^{p}$-norm over rings was given by Johnson \cite{Johnson} which is more extensively used. Currently, several versions of classical spaces, such as Sobolev, Hardy, Triebel-Lizorkin and Besov based on Herz spaces have already been introduced in the literature due to the growing interest in harmonic analysis and PDEs in those frameworks (see, e.g., \cite{Lucas_BWH, Cuerva_Herrero, Grafakos_Li_Yang, Tsutsui_2011, Xu_2005} and their references). Moreover, we point out that the standard Euler equations (constant density case) were treated in Besov-Herz spaces by \cite{Lucas_Herz}.

The nonhomogeneous Besov-Herz space $BK_{p,q,r}^{\alpha ,s}$ is larger than the nonhomogeneous Besov space $B_{p,r}^{s}$ for $s\in \mathbb{R}$, $\alpha=0 $, $q=\infty $, $1\leq p<\infty ,$ and $1\leq r\leq \infty $ (see Remark \ref{Rem:comp-spaces}). Then, for a fixed index $s$, we are able to construct a well-defined flow starting from a larger initial-data class for both the velocity and density. Even for the initial density $a_{0}$, the smallness condition is taken in the weaker norm of the Besov-Herz space which allows us to consider some large data in other setting such as Sobolev $H_{p}^{s},W^{s,p}$ and Besov $B_{p,q}^{s}$ spaces.

Below we state our main results.

\begin{theorem} \label{the:Euler}
Consider $1<p<\infty $, $1\leq q,r\leq \infty $, $0\leq\alpha < n(1-1/p),$ and also $s\geq n/p+1.$ Suppose also that $r=1$ when $s=n/p+1$.

\begin{itemize}
\item[(i)] Let $a_{0}\in BK_{p,q,r}^{\alpha ,s}$, $f\in L_{T}^{1}(BK_{p,q,r}^{\alpha ,s})$ and $u_{0}\in BK_{p,q,r}^{\alpha ,s}$ with $\mathrm{div}\hspace{0.05cm}u_{0}=0$. There exist $T\in (0,\infty )$ and a small constant $c>0$ such that if $\Vert a_{0}\Vert_{BK_{p,q,r}^{\alpha ,s}}\leq c,$ then system (\ref{sist:Euler_1}) admits a unique solution $(a,u,\nabla \pi )$ satisfying
\begin{equation} \label{aux-sol-1}
a\in C([0,T];BK_{p,q,r}^{\alpha ,s}),\hspace{0.3cm}u\in C([0,T];BK_{p,q,r}^{\alpha ,s})\hspace{0.3cm}\text{and}\hspace{0.3cm}\nabla\pi \in L_{T}^{1}(BK_{p,q,r}^{\alpha ,s}).
\end{equation}

\item[(ii)] Let $\{(a_{0,k},u_{0,k})\}_{k\in \mathbb{N}}$ be a bounded sequence of pairs in $BK_{p,q,r}^{\alpha ,s}$ such that $a_{0,k}\rightarrow a_{0}$ and $u_{0,k}\rightarrow u_{0}$ in $BK_{p,q,r}^{\alpha ,s-1}$ as $k\rightarrow \infty $. Consider $(a_{k},u_{k})$ and $(a,u)$ the solutions obtained in item (i) with the respective initial data $(a_{0,k},u_{0,k})$ and $(a_{0},u_{0})$. Then, there exists $T>0$ such that $\{(a_{k},u_{k})\}_{k\in \mathbb{N}}$ is bounded in $L_{T}^{\infty}(BK_{p,q,r}^{\alpha,s})\times L_{T}^{\infty }(BK_{p,q,r}^{\alpha ,s})$ and
\begin{equation}
(a_{k},u_{k})\overset{k\rightarrow \infty }{\longrightarrow }(a,u)\;\;\;\text{in}\;\;C([0,T];BK_{p,q,r}^{\alpha ,s-1})\times C([0,T];BK_{p,q,r}^{\alpha ,s-1}).
\end{equation}
\end{itemize}
\end{theorem}

\vspace{0.2cm}

\begin{remark} \label{Rem:blow-up}
(Blow-up criterion) Let $0<T^{\ast }<\infty $ and $a_{0},u_{0},f$ as in Theorem \ref{the:Euler}. Then, the corresponding local solution $(a,u,\nabla \pi )$ blows up at time $T^{\ast }>T$ in $BK_{p,q,r}^{\alpha ,s},$ that is,
\begin{align*}
\limsup_{t\nearrow T^{\ast }}\Vert a(t)\Vert _{BK_{p,q,r}^{\alpha ,s}}=\infty  \hspace{0.5cm} or \hspace{0.5cm} \limsup_{t\nearrow T^{\ast}}\Vert u(t)\Vert _{BK_{p,q,r}^{\alpha ,s}}=\infty,
\end{align*}
if and only if
\begin{equation*}
\int_{0}^{T^{\ast }}\Vert \nabla \times u(t)\Vert _{\dot{B}_{\infty ,\infty}^{0}}dt=\infty \hspace{0.35cm}\left( \text{resp.}\int_{0}^{T^{\ast }}\Vert\nabla \times u(t)\Vert _{\dot{B}_{\infty ,1}^{0}}dt=\infty \right) ,
\end{equation*}
when $s>n/p+1$ with $1\leq r\leq \infty $ (resp. $s=n/p+1$ with $r=1$). The proof follows essentially by proceeding as in the proof of the linear estimates in Propositions \ref{prop:transport} and \ref{prop:Euler} combined with a logarithmic inequality in the Besov-Herz setting, namely
\begin{align*}
\Vert u\Vert _{L^{\infty }}\lesssim 1+\Vert u\Vert_{\dot{B}_{\infty ,\infty}^{0}}\left( \log ^{+}\Vert u\Vert_{BK_{p,q,r}^{\alpha ,s^{\prime}}}+1\right)
\end{align*}
for some $s^{\prime }>n/p$, or estimate $\Vert \nabla u\Vert_{L^{\infty }}\lesssim \Vert \nabla u\Vert _{\dot{B}_{\infty ,1}^{0}}$, according to the corresponding case. Moreover, we need the decomposition $\nabla u=\mathcal{P}w+Mw$, where $w=\nabla \times u$, $\mathcal{P}$ is a singular integral operator and $M$ is a constant matrix.
\end{remark}

\begin{remark} \

\begin{itemize}
\item[(i)] Here we focus on the case $n\geq3$. However, with a slight adaptation of the proofs, Theorem \ref{the:Euler} holds for $n=2$.

\item[(ii)] Adapting some arguments in \cite{Danchin_2010} and \cite{Danchin-2011} to our framework, we believe that the smallness condition on $a_0$ in Theorem \ref{the:Euler} could be removed by considering an additional $L^2$-restriction on the potential part of $f$ and a $L^2$-setting for the pressure term $\nabla \pi$, as well as some further constraints on the indexes of the functional spaces.
\end{itemize}

\end{remark}

In what follows, comparing with previous works (see, e.g., \cite{Chae_2003_2, Danchin_2010, Danchin-2011, Wei_2013, Zhou_2010} and their references), we discuss some difficulties that we needed to overcome to achieve our results. In view of the presence of the base-space $K_{p,q}^{\alpha },$ a significant part of the $L^{p}(\mathbb{R}^{n})$-theory does not work in a simple way in the framework of Besov-Herz spaces. As a matter of fact, the $K_{p,q}^{\alpha }$-norm involves a localized handling of $L^{p}$-norms (in original variables $x$) by means of both ball $A_{-1}$ and rings $A_{k}$ of $\mathbb{R}^{n}$ (see Definition \ref{def:herz_spaces}). Thus, key ingredients employed in Besov spaces, such as integration by parts and energy-like arguments, are difficult to implement in our setting. In order to overcome these difficulties and develop the needed linear estimates (see Propositions \ref{prop:transport} and \ref{prop:Euler}), we handle both velocity and density equations by means of the volume-preserving map $X$ associated with $u$, inspired by some arguments that can be found in \cite{Chae_2003_3, Chae-2004, Chemin_2, Lucas_Herz, Takada, Vishik-1998}. Furthermore, we derive some commutator estimates useful for our ends, involving the velocity $u$, the pressure $\pi $ and the density $a$ (see Lemmas \ref{lem:comutator_uv} and \ref{lem:comutator_pressure}) by adapting some estimates of previous works (see, e.g., \cite{Danchin_2010, Lucas_Herz, Zhou_2010}). Also, we use some Bernstein-type inequalities on Besov-Herz and Sobolev-Herz spaces, mainly to estimate the pressure term $\nabla \pi $ (see Remark \ref{Rem:bernstein_bh}).

The outline of this manuscript is as follows. Section \ref{sec:02} is devoted to some preliminaries by introducing notations and reviewing definitions and basic properties on Herz, Sobolev-Herz and Besov-Herz spaces, as well as suitable estimates for volume-preserving maps. The subject of Section \ref{sec:com_prod} is a set of estimates for the product and commutator operators. In Section \ref{sec:linear_est} we derive core estimates for some linear conservation-law systems linked to structure of (\ref{sist:Euler_1}). In Section \ref{sec:existence}, by means of an approximate linear problem and uniform estimates for it, as well as some contraction estimates, we perform the proof of Theorem \ref{the:Euler} through six subsections.

\section{Preliminaries} \label{sec:02}

In this section we recall some definitions and basic properties about some operators and functional spaces that will be useful in the course of this manuscript.

The notation $A\lesssim B$ means that there is a generic constant $C>0$, which may be different on different lines, such that $A\leq CB$. Also, we use $A\approx B$ when $A\lesssim B$ and $B\lesssim A$. The commutator between two operators $F_{1}$ and $F_{2}$ is denoted by $[F_{1},F_{2}]=F_{1}F_{2}-F_{2}F_{1}$. Moreover, $C_{0}^{\infty }(\mathbb{R}^{n}),$ $\mathcal{S}=\mathcal{S}(\mathbb{R}^{n})$, and $\mathcal{S}^{\prime}=\mathcal{S}^{\prime }(\mathbb{R}^{n})$ stand for the space of compactly supported smooth functions, Schwartz space, and the space of tempered distributions, respectively. Also, consider $\mathbb{N}_{0}:= \mathbb{N} \cup \{0\}$ to represent the set of non-negative integers.

It is worth remembering that in 2005 H. Xu \cite{Xu_2005} introduced the spaces $\dot{B}K_{p,q,r}^{\alpha,s}$ and $BK_{p,q,r}^{\alpha,s}$ called homogeneous and nonhomogeneous Besov-Herz spaces, which are Besov spaces based on Herz $K_{p,q}^{\alpha}$ spaces. In the sequel we give some preliminaries on Herz spaces. For more details, see \cite{Lucas_Herz, Grafakos_Li_Yang, Herz, Johnson, Li_Yang, Xu_2005}.

\begin{definition} \label{def:herz_spaces}
Let $\alpha \in \mathbb{R}$ and $1\leq p,q\leq\infty $. The nonhomogeneous Herz space $K_{p,q}^{\alpha }=K_{p,q}^{\alpha }(\mathbb{R}^{n})$ is defined as the set of all functions $u\in L_{\mathrm{loc}}^{p}(\mathbb{R}^{n})$ such that
\begin{equation*}
\Vert u\Vert _{K_{p,q}^{\alpha }}:=\left( \sum_{k\geq -1}\left( 2^{\alpha k}\Vert u\Vert _{L^{p}(A_{k})}\right) ^{q}\right) ^{1/q}<\infty ,
\end{equation*}
with the usual modification for $p=\infty $ or $q=\infty $, where
\begin{equation*}
A_{-1}:=\{x\in \mathbb{R}^{n}\;:\;|x|<2^{-1}\} \hspace{0.5cm}\text{and} \hspace{0.5cm} A_{k}:=\{x\in \mathbb{R}^{n}\;:\;2^{k-1}\leq |x|<2^{k}\}, \hspace{0.25cm} \text{for}\hspace{0.25cm} k\geq 0.
\end{equation*}
\end{definition}

\begin{remark}
It is not difficult to prove that $K_{p,q}^{\alpha}$ is a Banach space with the norm $\Vert\cdot\Vert_{K_{p,q}^{\alpha}}$ for $1\leq p<\infty$, $1\leq q\leq\infty$ and $\alpha\in\mathbb{R}$, see e.g. \cite{Hernandez}. Moreover, we consider $K_{1,q}^{\alpha}$ as a space of signed measures where $\Vert u \Vert_{L^1(A_k)}$ denote the total variation of $u$ on $A_k$.
\end{remark}

In what follows, we recall the Littlewood-Paley decompositions, see e.g. \cite{Chemin_2, Lemarie} for further details. For that, consider $\psi $ and $\varphi \in C_{0}^{\infty }(\mathbb{R}^{n})$ radially symmetric functions supported in the closed ball $\mathcal{B}:=\{\xi \in \mathbb{R}^{n}:|\xi|\leq 4/3\}$ and in the ring $\mathcal{C}:=\{\xi \in \mathbb{R}^{n}:3/4\leq|\xi |\leq 8/3\},$ respectively, such that
\begin{equation} \label{decomp:unitaria}
\sum_{j\in \mathbb{Z}}\varphi _{j}(\xi )=1,\hspace{0.25cm}\text{for }\xi \in\mathbb{R}^{n}/\{0\},\text{ \ \ \ and \ \ \ }\psi (\xi )+\sum_{j\in \mathbb{N}_{0}}\varphi _{j}(\xi )=1,\text{ for }\xi \in \mathbb{R}^{n},
\end{equation}
where $\varphi _{j}(\xi ):=\varphi (2^{-j}\xi )$. Now, for $u\in \mathcal{S}^{\prime }(\mathbb{R}^{n})$, we define the localization operator $\dot{\Delta}_{j}$ and $\Delta _{j}u$ by
\begin{equation*}
\dot{\Delta}_{j}u:=\mathcal{F}^{-1}[\varphi _{j}\mathcal{F}u],\hspace{0.25cm} \mbox{for all $j\in\mathbb{Z}$} \hspace{0.5cm}\mbox{and}\hspace{0.5cm}\Delta_{j}u := \left\{
\begin{array}{cl}
0, & j\leq -2, \\
\mathcal{F}^{-1}[\psi \mathcal{F}u], & j=-1, \\
\dot{\Delta}_{j}u, & j\geq 0.
\end{array} \right.
\end{equation*}
Also, consider the low-frequency cut-off $S_{j}u:=\sum_{l\leq j-1}\Delta _{l}u$, where $\widehat{u}=\mathcal{F}u$ denotes the Fourier transform of $u$ on $\mathbb{R}^{n}$ and $\mathcal{F}^{-1}u$ denotes the inverse Fourier transform. Then, we have the Littlewood-Paley decompositions
\begin{equation*}
u=\sum_{j\in \mathbb{Z}}\dot{\Delta}_{j}u,\hspace{0.25cm} \mbox{for all $u\in\mathcal{S}'(\mathbb{R}^n)/\mathcal{P}[\mathbb{R}^n]$,} \hspace{0.75cm} \mbox{and}\hspace{0.75cm}u=\sum_{j\in \mathbb{Z}}\Delta _{j}u,\hspace{0.25cm} \mbox{for all $u \in \mathcal{S}'(\mathbb{R}^n)$},
\end{equation*}
where $\mathcal{P}=\mathcal{P}[\mathbb{R}^{n}]$ denotes the set of polynomials with $n$ variables. Moreover, the above decompositions satisfy
\begin{equation*}
\Delta _{j}(\Delta _{k}u)\equiv 0,\hspace{0.25cm}\mbox{if $|j-k|\geq2$,} \hspace{0.5cm} \mbox{and}\hspace{0.5cm}\Delta _{j}(S_{k-1}u\Delta _{k}u)\equiv 0,\hspace{0.25cm}\mbox{if $|j-k|\geq5$}.
\end{equation*}

With the above definitions in hand, we can recall homogeneous and nonhomogeneous Besov-Herz spaces (see, e.g., \cite{Lucas_Herz, Xu_2005}).

\begin{definition} \label{def:Besov_Herz}
Let $1\leq p,q,r\leq \infty $, and $\alpha ,s\in\mathbb{R}$. The homogeneous Besov-Herz space $\dot{B}K_{p,q,r}^{\alpha ,s} = \dot{B}K_{p,q,r}^{\alpha ,s}(\mathbb{R}^{n})$ is defined as the set of all $u\in \mathcal{S}^{\prime }/\mathcal{P}$ such that $\dot{\Delta}_{j}u\in K_{p,q}^{\alpha }$ for every $j\in \mathbb{Z}$, and that
\begin{equation*}
\Vert u\Vert _{\dot{B}K_{p,q,r}^{\alpha ,s}}:=\left\{
\begin{array}{l}
\left( \displaystyle\sum_{j\in \mathbb{Z}}(2^{sj}\Vert \dot{\Delta}_{j}u\Vert _{K_{p,q}^{\alpha }})^{r}\right) ^{1/r}<\infty ,\hspace{0.5cm} \mbox{if $r<\infty$}, \\
\sup_{j\in \mathbb{Z}}2^{sj}\Vert \dot{\Delta}_{j}u\Vert _{K_{p,q}^{\alpha}}<\infty ,\hspace{0.5cm}\mbox{if $r=\infty$}.%
\end{array} \right.
\end{equation*}
The nonhomogeneous Besov-Herz space $BK_{p,q,r}^{\alpha,s}=BK_{p,q,r}^{\alpha ,s}(\mathbb{R}^{n})$ is defined as the set of all $u\in \mathcal{S}^{\prime }$ such that $\Delta _{j}u\in K_{p,q}^{\alpha}$ for every $j\geq -1$, and that
\begin{equation*}
\Vert u\Vert _{BK_{p,q,r}^{\alpha ,s}}:=\left\{
\begin{array}{l}
\left( \displaystyle\sum_{j\geq -1}(2^{sj}\Vert \Delta_{j}u \Vert_{K_{p,q}^{\alpha }})^{r}\right) ^{1/r}<\infty ,\hspace{0.5cm} \mbox{if $r<\infty$}, \\
\sup_{j\geq -1}2^{sj}\Vert \Delta _{j}u\Vert _{K_{p,q}^{\alpha }}<\infty , \hspace{0.5cm}\mbox{if $r=\infty$}.
\end{array} \right.
\end{equation*}
\end{definition}

In the sequel we recall the definition of homogeneous Sobolev-Herz spaces. For similar definitions based on $L^{p}$-spaces, weak-Herz spaces, or nonhomogeneous Sobolev-Herz spaces, see e.g. \cite{Bergh_Lofstrom, Lucas_BWH, Lucas_Herz, Xu_Yang_2003}.

\begin{definition} \label{def:Sobolev_Herz}
Let $1\leq p,q,r\leq \infty $, and $\alpha ,s\in\mathbb{R}$. The homogeneous Sobolev-Herz space $\dot{H}K_{p,q}^{\alpha ,s} = \dot{H}K_{p,q}^{\alpha ,s}(\mathbb{R}^{n})$ is the set of all $u\in \mathcal{S}^{\prime }/\mathcal{P}$ such that
\begin{equation*}
\Vert u\Vert _{\dot{H}K_{p,q}^{\alpha ,s}}:=\Vert I^{s}u \Vert_{K_{p,q}^{\alpha }}<\infty ,
\end{equation*}
where $I^{s}u=\mathcal{F}^{-1}[|\cdot |^{s}\mathcal{F}u]$ with $I^{-s}$ being the Riesz potential of order $s$. In particular, if $s=0$ we have $\dot{H}K_{p,q}^{\alpha ,0} \equiv K_{p,q}^{\alpha}$ with equivalent norms.
\end{definition}

Next, we highlight some properties about Sobolev and Besov-Herz spaces. For $p<\infty $, the spaces $\dot{B}K_{p,q,r}^{\alpha ,s}$, $BK_{p,q,r}^{\alpha,s}$ and $\dot{H}K_{p,q}^{\alpha ,s}$ are Banach spaces and satisfy the inclusions $\dot{B}K_{p,q,r}^{\alpha ,s},\dot{H}K_{p,q}^{\alpha ,s}\subset\mathcal{S}^{\prime }/\mathcal{P}$ and $BK_{p,q,r}^{\alpha ,s}\subset
\mathcal{S}^{\prime }.$ Moreover, for $s>0$ we have that
\begin{equation} \label{def:new_norm_inhomog}
BK_{p,q,r}^{\alpha ,s}=\dot{B}K_{p,q,r}^{\alpha ,s}\cap K_{p,q}^{\alpha } \hspace{0.5cm} \text{with} \hspace{0.5cm} \Vert u \Vert_{BK_{p,q,r}^{\alpha,s}}\approx \Vert u\Vert _{\dot{B}K_{p,q,r}^{\alpha ,s}}+\Vert u\Vert_{K_{p,q}^{\alpha }}.
\end{equation}

Using the decompositions given in (\ref{decomp:unitaria}), it is possible to show the following embeddings. Let $\alpha \geq 0$ and $1\leq p,q\leq \infty $. For $s>n/p$ and $1\leq r\leq \infty $, or $s=n/p$ and $r=1$, it follows that
\begin{equation} \label{l_infty_imersion}
BK_{p,q,r}^{\alpha ,s}\hookrightarrow L^{\infty }\hspace{0.5cm}\text{and} \hspace{0.5cm} \dot{B}K_{p,q,1}^{\alpha ,n/p}\hookrightarrow L^{\infty}.
\end{equation}%
Moreover, if $s_{1}>s_{2}$ we have the inclusion relation $BK_{p,q,r_{1}}^{\alpha ,s_{1}}\subset BK_{p,q,r_{2}}^{\alpha,s_{2}}$, for all $1\leq r_{1},r_{2}\leq \infty $.

In the next remark we observe some relations between Herz $K_{p,q}^{\alpha }$, Lebesgue $L^{p}$, Morrey $M_{p,\lambda }$ and Besov $B_{p,r}^{s}$ spaces (see \cite[Lemma 3.3]{Lucas_Herz}). For definitions and basic properties about the spaces $B_{p,r}^{s}$ and $M_{p,\lambda },$ the reader is referred to \cite{Bergh_Lofstrom} and \cite{Kato}, respectively.

\begin{remark} \label{Rem:comp-spaces}
For $1\leq p<\infty ,$ $1\leq r\leq \infty ,$ and $s\in \mathbb{R},$ we have the strict inclusions $L^{p}\subset K_{p,\infty}^{0}$ and $B_{p,r}^{s}\subset BK_{p,\infty ,r}^{0,s}$. Moreover, if $1\leq q<\infty $ and $0\leq \lambda <n$ with $\lambda \neq n(1-q/p)$ when $q<p$, then $K_{p,\infty }^{0}\not\subset M_{q,\lambda }$.
\end{remark}

To conclude this section, we present the decomposition of the Bony paraproduct, see more details in \cite{Bony}. Let $u,v\in \mathcal{S}^{\prime }$. The product of $u$ by $v$ can be written as
\begin{equation} \label{decomp_para_product}
uv=T_{u}v+T_{v}u+\mathcal{R}(u,v)=T_{u}v+R(u,v),
\end{equation}
where
\begin{equation*}
\begin{array}{rclcrcl}
\displaystyle T_{u}v & := & \displaystyle\sum_{j\in \mathbb{Z}}S_{j-1}u\Delta _{j}v, & \hspace{1cm} & R(u,v) & := & \displaystyle \sum_{j\in \mathbb{Z}}\Delta _{j}uS_{j+2}v, \\
\mathcal{R}(u,v) & := & \displaystyle\sum_{j\in \mathbb{Z}}\Delta _{j}u\widetilde{\Delta }_{j}v, & \hspace{1cm} & \widetilde{\Delta }_{j}v & := & \displaystyle\sum_{i=-1}^{1}\Delta _{j-i}v.
\end{array}
\end{equation*}

\subsection{Some basic estimates in Herz spaces}

We start by recalling H\"{o}lder-type inequality in the framework of Herz spaces (see, e.g., \cite{Tsutsui_2011}).

\begin{lemma} \label{lem:holder_herz}
Let $1\leq p,p_{1},p_{2},q,q_{1},q_{2}\leq \infty $ and $\alpha ,\alpha _{1},\alpha _{2}\in \mathbb{R}$ satisfy $1/p=1/p_{1}+1/p_{2}$, $1/q=1/q_{1}+1/q_{2}$ and $\alpha =\alpha _{1}+\alpha_{2}$. Then, there is a constant $C>0$ such that
\begin{equation*}
\Vert uv\Vert _{K_{p,q}^{\alpha }}\leq C\Vert u \Vert_{K_{p_{1},q_{1}}^{\alpha _{1}}} \Vert v \Vert_{K_{p_{2},q_{2}}^{\alpha_{2}}},
\end{equation*}
for all $u\in K_{p_{1},q_{1}}^{\alpha _{1}}$ and $v\in K_{p_{2},q_{2}}^{\alpha _{2}}$. In particular, if $u\in L^{\infty }$ and $v\in K_{p,q}^{\alpha }$, we have that
\begin{equation} \label{est:holder_herz}
\Vert uv\Vert _{K_{p,q}^{\alpha }}\leq C\Vert u\Vert_{L^{\infty}}\Vert v\Vert _{K_{p,q}^{\alpha }}.
\end{equation}
\end{lemma}

Under suitable conditions, we also have a Young-type inequality in Herz spaces (see, e.g., \cite{Lucas_Herz}).

\begin{lemma}
\label{lem:young_herz} Let $\alpha \in \mathbb{R}$ and $1\leq p,q\leq \infty
$. Consider $\phi \in L^{1}$ satisfying%
\begin{equation*}
M_{\phi }:=\left\{
\begin{array}{l}
\max \left\{ \Vert \phi \Vert _{L^{1}},\Vert |\cdot |^{\beta }\phi \Vert
_{L^{1}},\Vert |\cdot |^{2\beta +\alpha p}\phi \Vert _{L^{1}}\right\}
<\infty ,\;\;\text{if}\;\alpha \geq 0, \\
\max \left\{ \Vert \phi \Vert _{L^{1}},\Vert |\cdot |^{\beta -\alpha p}\phi
\Vert _{L^{1}},\Vert |\cdot |^{2\beta }\phi \Vert _{L^{1}}\right\} <\infty
,\;\;\text{if}\;\alpha <0,%
\end{array}%
\right.
\end{equation*}%
for some $\beta >0$. Then, there is a constant $C>0$, independent of $%
\varphi $, such that
\begin{equation}
\Vert \phi \ast u\Vert _{K_{p,q}^{\alpha }}\leq CM_{\phi }\Vert u\Vert
_{K_{p,q}^{\alpha }},  \label{est-young1}
\end{equation}%
for all $u\in K_{p,q}^{\alpha }$.
\end{lemma}

The following lemma compiles some Bernstein-type inequalities in the
framework of Herz spaces found in the literature (see {\cite{Lucas_Herz}}).
For $0<R_{1}<R_{2},$ consider the notations
\begin{align*}
\mathcal{B}(0,R_{1}):=\left\{\xi \in \mathbb{R}^{n};\left\vert \xi \right\vert \leq R_{1}\right\} \hspace{0.5cm} \mbox{and} \hspace{0.5cm} \mathcal{C}(0,R_{1},R_{2}):=\left\{ \xi \in \mathbb{R}^{n};R_{1}\leq \left\vert \xi \right\vert \leq R_{2}\right\}.
\end{align*}

\begin{lemma}
\label{lem:bernstein_herz} Let $0<R_{1}<R_{2},$ $1\leq p,q\leq \infty $, $%
\alpha \in \mathbb{R}$, and $j\in \mathbb{Z}$.

\begin{itemize}
\item[(i)] For $u\in K_{p,q}^{\alpha }$ satisfying $\mathrm{supp}\hspace{%
0.05cm}\widehat{u}\subset \mathcal{B}(0,R_{1}2^{j})$, it follows that
\begin{equation}
\Vert \partial ^{\beta }u\Vert _{K_{p,q}^{\alpha }}\leq C2^{j|\beta |}\Vert
u\Vert _{K_{p,q}^{\alpha }},  \label{est:bernstein_herz_B}
\end{equation}%
for some constant $C>0$ independent of $j$.

\item[(ii)] For $u\in K_{p,q}^{\alpha }$ satisfying $\mathrm{supp}\hspace{%
0.05cm}\widehat{u}\subset \mathcal{C}(0,R_{1}2^{j},R_{2}2^{j})$, it follows
that
\begin{equation*}
C^{-1}2^{j|\beta |}\Vert u\Vert _{K_{p,q}^{\alpha }}\leq \Vert \partial
^{\beta }u\Vert _{K_{p,q}^{\alpha }}\leq C2^{j|\beta |}\Vert u\Vert
_{K_{p,q}^{\alpha }},
\end{equation*}%
for some constant $C>0$ independent of $j$.

\item[(iii)] For $\alpha \geq 0$ and $u\in K_{p,q}^{\alpha }$ satisfying $%
\mathrm{supp}\hspace{0.05cm}\widehat{u}\subset \mathcal{B}(0,R_{1}2^{j})$,
it follows that
\begin{equation}
\Vert u\Vert _{L^{\infty }(\mathbb{R}^{n})}\leq C2^{jn/p}\Vert u\Vert
_{K_{p,q}^{\alpha }},  \label{est:l_infty_herz}
\end{equation}%
for some constant $C>0$ independent of $j$.
\end{itemize}
\end{lemma}

As a direct consequence of Lemma \ref{lem:bernstein_herz}, we have the
following remark.

\begin{remark}
\label{Rem:bernstein_bh} For $k\in \mathbb{N}_{0}$ and $\alpha ,s\in \mathbb{%
R}$, we have the estimates
\begin{align}
\hspace{1cm}\Vert \partial ^{\beta }u\Vert _{BK_{p,q,r}^{\alpha ,s}}& \leq
C\Vert u\Vert _{BK_{p,q,r}^{\alpha ,s+k}},\hspace{0.5cm}%
\mbox{for\;
$|\beta|=k$}  \label{bernstein_bh_inhomog} \\
C^{-1}\Vert u\Vert _{\dot{B}K_{p,q,r}^{\alpha ,s+k}}\leq \Vert \partial
^{\beta }u& \Vert _{\dot{B}K_{p,q,r}^{\alpha ,s}}\leq C\Vert u\Vert _{\dot{B}%
K_{p,q,r}^{\alpha ,s+k}},\hspace{0.5cm}\mbox{for\; $|\beta|=k$}
\label{bernstein_bh_homog}
\end{align}%
where $1\leq p,q\leq \infty $ and $C>0$ is a universal constant.
\end{remark}

Moreover, we have estimates similar to (\ref{bernstein_bh_homog}) in the
context of Sobolev-Herz spaces, namely
\begin{equation}
C^{-1}\Vert u\Vert _{\dot{H}K_{p,q}^{\alpha ,s+k}}\leq \sup_{|\beta
|=k}\Vert \partial ^{\beta }u\Vert _{\dot{H}K_{p,q}^{\alpha ,s}}\leq C\Vert
u\Vert _{\dot{H}K_{p,q}^{\alpha ,s+k}},  \label{bernstein_sh_homog}
\end{equation}%
where $1<p<\infty ,$ $1\leq q\leq \infty $, $k\in \mathbb{N}_{0}$, $s\in
\mathbb{R}$ and $0\leq \alpha <n(1-1/p).$ The proof of (\ref%
{bernstein_sh_homog}) can be obtained by following arguments in Sobolev $%
\dot{H}_{p}^{s}$ and Triebel-Lizorkin $\dot{F}_{p,q}^{s}$ spaces, as well as
equivalence relations between them, and some basic properties of Fourier
multipliers. For further details, the reader is referred to \cite{Grafakos,
Triebel, Xu_Yang_2003, Xu_Yang_2005}.

\subsection{Volume-preserving maps}

This part is devoted to recall an estimate for volume-preser- ving maps $X$ in
Herz spaces. For the proof, we refer the reader to \cite{Lucas_Herz}.

\begin{lemma}
\label{lem:field_x} Let $1\leq p,q\leq \infty $ and $\alpha \geq 0$. Assume
that $X:\mathbb{R}^{n}\rightarrow \mathbb{R}^{n}$ is a volume-preserving
diffeomorphism satisfying
\begin{equation}
|X^{\pm 1}(x_{0})-x_{0}|\leq \gamma ,\hspace{0.5cm}\text{for all }x_{0}\in
\mathbb{R}^{n},  \label{est:field_x}
\end{equation}%
and some fixed $\gamma >0$. Then, there exists a constant $C>0$ such that
\begin{equation}
C^{-1}\Vert u\Vert _{K_{p,q}^{\alpha }}\leq \Vert u\circ X\Vert
_{K_{p,q}^{\alpha }}\leq C\Vert u\Vert _{K_{p,q}^{\alpha }},
\label{est:field_composition}
\end{equation}%
for all $u\in K_{p,q}^{\alpha }$.
\end{lemma}

For the sake of handling, the remark below will be helpful in our context.

\begin{remark}
\label{rem:field_X} If $X$ is the flow generated by a field $v$ associated
with the ODE
\begin{equation}
\left\{
\begin{array}{l}
\partial _{t}X(y,t)=v(X(y,t),t) \\
X(y,0)=y%
\end{array}%
,\right. \;\;\;\;(y,t)\in \mathbb{R}^{n}\times \mathbb{R}^{+},
\label{sist:field_X}
\end{equation}%
then, $X$ verifies the following estimate
\begin{equation*}
|X^{\pm 1}(y,t)-y|\leq \int_{0}^{t}\Vert v(\tau )\Vert _{L^{\infty }}\;d\tau
.
\end{equation*}%
Moreover, since $\mathrm{div}\;v=0$, it follows that $X$ is a
volume-preserving diffeomorphism for each $t\geq 0$ and, consequently,
satisfies estimate (\ref{est:field_x}) under suitable time-integrability
conditions on $v.$
\end{remark}

\section{Commutator and product estimates} \label{sec:com_prod}

The present section contains estimates for product and commutator operators
in the context of Herz and Besov-Herz spaces, see e.g. \cite{Danchin_2010,
Lucas_Herz, Zhou_2010} for similar estimates.

We start with estimates for the product in Besov-Herz spaces.

\begin{lemma}
\label{lem:product_estimates} Let $1\leq p,q,r\leq \infty $ and $\alpha \in
\mathbb{R}$.

\begin{itemize}
\item[(i)] For $s>0,$ the following estimate
\begin{equation}
\Vert uv\Vert _{BK_{p,q,r}^{\alpha ,s}}\leq C\left( \Vert u\Vert _{L^{\infty
}}\Vert v\Vert _{BK_{p,q,r}^{\alpha ,s}}+\Vert v\Vert _{L^{\infty }}\Vert
u\Vert _{BK_{p,q,r}^{\alpha ,s}}\right)  \label{est:prod_uv_1}
\end{equation}%
holds true, for all $u,v\in BK_{p,q,r}^{\alpha ,s}\cap L^{\infty },$ where $%
C>0$ is a universal constant. The same estimate still holds in homogeneous
Besov-Herz spaces.

\item[(ii)] For $p<\infty $, $\alpha \geq 0$ and $s\geq n/p+1$, with $r=1$
if $s=n/p+1$, we have that
\begin{equation}
\Vert u\cdot \nabla v\Vert _{BK_{p,q,r}^{\alpha ,s-1}}\leq C\Vert u\Vert
_{BK_{p,q,r}^{\alpha ,s-1}}\Vert v\Vert _{BK_{p,q,r}^{\alpha ,s}},
\label{est:prod_uv_2}
\end{equation}%
for all $u\in BK_{p,q,r}^{\alpha ,s-1}$ and $v\in BK_{p,q,r}^{\alpha ,s},$
where $C>0$ is a universal constant. Moreover, assuming further {$\mathrm{div%
}\hspace{0.05cm}v=0,$ }the estimate
\begin{equation}
\Vert \mathrm{div}(u\cdot \nabla v)\Vert _{BK_{p,q,r}^{\alpha ,s-1}}\leq
C\Vert u\Vert _{BK_{p,q,r}^{\alpha ,s}}\Vert v\Vert _{BK_{p,q,r}^{\alpha ,s}}
\label{est:prod_uv_3}
\end{equation}%
holds true, for all $u,v\in BK_{p,q,r}^{\alpha ,s},$ where $C>0$ is a
universal constant.
\end{itemize}
\end{lemma}

\begin{remark}
The proof of Lemma \ref{lem:product_estimates} can be found in \cite%
{Lucas_Herz}, except for the case in homogeneous spaces, which can be proved
directly by using the Bony paraproduct decomposition (\ref%
{decomp_para_product}). Furthermore, observing that (\ref{est:prod_uv_1})
holds for $s-1$, $\mathrm{div}\hspace{0.05cm}v=0$, and
\begin{equation*}
\mathrm{div}(u\cdot \nabla v)=\nabla u:\nabla v=\sum_{i,j=1}^{n}\partial
_{i}u_{j}\partial _{i}v_{j},
\end{equation*}%
we obtain (\ref{est:prod_uv_3}) naturally from the embedding $%
BK_{p,q,r}^{\alpha ,s}\hookrightarrow L^{\infty }$ given in (\ref%
{l_infty_imersion}).
\end{remark}

Next, we present commutator estimates that can be seen as adaptations of
some previous estimates found in the literature, see \cite{Danchin_2010,
Lucas_Herz, Zhou_2010} and references therein.

\begin{lemma}
\label{lem:comutator_uv} Let $1\leq p<\infty $, $1\leq q,r\leq \infty $, $%
\alpha \geq 0$ and let $v$ be a divergence-free vector field.

\begin{itemize}
\item[(i)] For $s>0,$ there exists a constant $C>0$ such that
\begin{equation}
\left( \sum_{j\geq -1}\left( 2^{sj}\Vert \lbrack \Delta _{j},v\cdot \nabla
]u\Vert _{K_{p,q}^{\alpha }}\right) ^{r}\right) ^{1/r}\leq C\left( \Vert
\nabla u\Vert _{L^{\infty }}\Vert v\Vert _{BK_{p,q,r}^{\alpha ,s}}+\Vert
\nabla v\Vert _{L^{\infty }}\Vert u\Vert _{BK_{p,q,r}^{\alpha ,s}}\right) ,
\label{est:comutator_uv_1}
\end{equation}
for all $u,v\in BK_{p,q,r}^{\alpha ,s}$ with $\nabla u,\nabla v\in L^{\infty
}.$

\item[(ii)] For $s>n/p+1$ or $s=n/p+1$ with $r=1$, we have the estimate
\begin{equation}
\left( \sum_{j\geq -1}\left( 2^{sj}\Vert \lbrack \Delta _{j},v\cdot \nabla
]u\Vert _{K_{p,q}^{\alpha }}\right) ^{r}\right) ^{1/r}\leq C\Vert u\Vert
_{BK_{p,q,r}^{\alpha ,s}}\Vert v\Vert _{BK_{p,q,r}^{\alpha ,s}},
\label{est:comutator_uv_2}
\end{equation}%
for all $u,v\in BK_{p,q,r}^{\alpha ,s},$ where $C>0$ is a universal
constant. Moreover, there exists a constant $C>0$ such that
\begin{equation}
\left( \sum_{j\geq -1}\left( 2^{(s-1)j}\Vert \lbrack \Delta _{j},v\cdot
\nabla ]u\Vert _{K_{p,q}^{\alpha }}\right) ^{r}\right) ^{1/r}\leq C\Vert
u\Vert _{BK_{p,q,r}^{\alpha ,s-1}}\Vert v\Vert _{BK_{p,q,r}^{\alpha ,s}},
\label{est:comutator_uv_3}
\end{equation}%
for all $u\in BK_{p,q,r}^{\alpha ,s-1}$ and $v\in BK_{p,q,r}^{\alpha ,s}$.
\end{itemize}
\end{lemma}

\textbf{Proof:} By the Bony decomposition (\ref{decomp_para_product}), we
can write
\begin{eqnarray}
\left[ \Delta _{j},v\cdot \nabla \right] u &=&\Delta _{j}(\mathcal{R}%
(v,\nabla u))+\Delta _{j}\left( T_{\nabla u}v\right) -R\left( v,\Delta
_{j}\nabla u\right) -\left[ T_{v},\Delta _{j}\right] \nabla u  \notag \\
&=:&\mathcal{R}_{j}^{1}+\mathcal{R}_{j}^{2}+\mathcal{R}_{j}^{3}+\mathcal{R}%
_{j}^{4}.  \label{est:comt_decomp}
\end{eqnarray}%
For $\mathcal{R}_{j}^{1}$, first note that $\mathcal{R}_{j}^{1}=\sum_{j-k%
\leq 3}\Delta _{j}(\Delta _{k}v\tilde{\Delta}_{k}\nabla u)$. \ By Lemma \ref%
{lem:young_herz} and H\"{o}lder inequality (\ref{est:holder_herz}), we have
\begin{equation*}
\Vert \mathcal{R}_{j}^{1}\Vert _{K_{p,q}^{\alpha }}\lesssim \sum_{j-k\leq
3}\Vert \Delta _{k}v\tilde{\Delta}_{k}\nabla u\Vert _{K_{p,q}^{\alpha
}}\lesssim \Vert \nabla u\Vert _{L^{\infty }}\sum_{j-k\leq 3}\Vert \Delta
_{k}v\Vert _{K_{p,q}^{\alpha }}.
\end{equation*}%
Consequently, for $m=j-k$, multiplying both sides by $2^{sj}$ and taking the
$\ell ^{r}$-norm, we arrive at
\begin{align}
\left( \sum_{j\geq -1}\left( 2^{sj}\Vert \mathcal{R}_{j}^{1}\Vert
_{K_{p,q}^{\alpha }}\right) ^{r}\right) ^{1/r}& \lesssim \Vert \nabla u\Vert
_{L^{\infty }}\left( \sum_{j\geq -1}\left( \sum_{m\leq 3}2^{sm}\left(
2^{s(j-m)}\Vert \Delta _{j-m}v\Vert _{K_{p,q}^{\alpha }}\right) \right)
^{r}\right) ^{1/r}  \notag \\
& \lesssim \Vert \nabla u\Vert _{L^{\infty }}\sum_{m\leq 3}2^{sm}\left(
\sum_{j\geq -1}\left( 2^{s(j-m)}\Vert \Delta _{j-m}v\Vert _{K_{p,q}^{\alpha
}}\right) ^{r}\right) ^{1/r},  \label{est:metade}
\end{align}%
where the last estimate follows from Minkowski inequality. Then, for $s>0$,
it follows that
\begin{equation}
\left( \sum_{j\geq -1}\left( 2^{sj}\Vert \mathcal{R}_{j}^{1}\Vert
_{K_{p,q}^{\alpha }}\right) ^{r}\right) ^{1/r}\lesssim \Vert \nabla u\Vert
_{L^{\infty }}\Vert v\Vert _{BK_{p,q,r}^{\alpha ,s}}.  \label{est:comt_R_1}
\end{equation}%
For $\mathcal{R}_{j}^{2},$ in view of the support of $\varphi _{j},$ we can
express $\mathcal{R}_{j}^{2}=\sum_{|j-k|\leq 4}\Delta _{j}\left(
S_{k-1}\nabla u\Delta _{k}v\right) $. Now, using again Lemma \ref%
{lem:young_herz} and H\"{o}lder inequality (\ref{est:holder_herz}), we can
estimate
\begin{equation*}
\Vert \mathcal{R}_{j}^{2}\Vert _{K_{p,q}^{\alpha }}\leq \sum_{|j-k|\leq
4}\Vert \Delta _{j}\left( S_{k-1}\nabla u\Delta _{k}v\right) \Vert
_{K_{p,q}^{\alpha }}\lesssim \Vert \nabla u\Vert _{L^{\infty
}}\sum_{|j-k|\leq 4}\Vert \Delta _{k}v\Vert _{K_{p,q}^{\alpha }}.
\end{equation*}%
Therefore, analogously as in (\ref{est:metade}), for $s\in \mathbb{R}$, it
follows that
\begin{equation}
\left( \sum_{j\geq -1}\left( 2^{sj}\Vert \mathcal{R}_{j}^{2}\Vert
_{K_{p,q}^{\alpha }}\right) ^{r}\right) ^{1/r}\lesssim \Vert \nabla u\Vert
_{L^{\infty }}\Vert v\Vert _{BK_{p,q,r}^{\alpha ,s}}.  \label{est:comt_R_2}
\end{equation}%
For $\mathcal{R}_{j}^{3}$, thanks to the support of $\mathcal{F}%
(S_{k+2}\Delta _{j}\nabla u)$, we can write $\mathcal{R}_{j}^{3}=-\sum_{j-k%
\leq 2}\Delta _{k}vS_{k+2}(\Delta _{j}\nabla u)$ and estimate it as
\begin{equation*}
\Vert \mathcal{R}_{j}^{3}\Vert _{K_{p,q}^{\alpha }}\lesssim \sum_{j-k\leq
2}\Vert \Delta _{k}vS_{k+2}(\Delta _{j}\nabla u)\Vert _{K_{p,q}^{\alpha
}}\lesssim \Vert \nabla u\Vert _{L^{\infty }}\sum_{j-k\leq 2}\Vert \Delta
_{k}v\Vert _{K_{p,q}^{\alpha }},
\end{equation*}%
which leads us to
\begin{equation}
\left( \sum_{j\geq -1}\left( 2^{sj}\Vert \mathcal{R}_{j}^{3}\Vert
_{K_{p,q}^{\alpha }}\right) ^{r}\right) ^{1/r}\lesssim \Vert \nabla u\Vert
_{L^{\infty }}\Vert v\Vert _{BK_{p,q,r}^{\alpha ,s}}.  \label{est:comt_R_3}
\end{equation}%
For the last term in (\ref{est:comt_decomp}), by changing variables and
convolution properties, we can write
\begin{equation*}
\lbrack T_{v},\Delta _{j}]\nabla u=\sum_{|j-k|\leq 4}2^{-j}\int_{\mathbb{R}%
^{n}}\varphi ^{\vee }(y)\int_{0}^{1}(y\cdot \nabla )S_{k-1}v(x-2^{-j}y\tau
)\;d\tau \;\Delta _{k}\nabla u(x-2^{-j}y)dy.
\end{equation*}%
Then, applying the Herz norm and using Young inequality in Lemma \ref%
{lem:young_herz}, since $\varphi \in \mathcal{S}$, we obtain that
\begin{equation*}
\Vert \mathcal{R}_{j}^{4}\Vert _{K_{p,q}^{\alpha }}\lesssim \sum_{|j-k|\leq
4}2^{-j}\Vert \nabla S_{k-1}v\Vert _{L^{\infty }}\Vert y\cdot \varphi ^{\vee
}\Vert _{L^{1}}\Vert \Delta _{k}\nabla u\Vert _{K_{p,q}^{\alpha }}\lesssim
\Vert \nabla v\Vert _{L^{\infty }}\sum_{|j-k|\leq 4}2^{k-j}\Vert \Delta
_{k}u\Vert _{K_{p,q}^{\alpha }},
\end{equation*}%
where above we used the Bernstein inequality (\ref{est:bernstein_herz_B}).
Then, making $m=j-k$ and taking $\ell ^{r}$-norm, for $s\in \mathbb{R}$, we
get
\begin{align}
\left( \sum_{j\geq -1}\left( 2^{sj}\Vert \mathcal{R}_{j}^{4}\Vert
_{K_{p,q}^{\alpha }}\right) ^{r}\right) ^{1/r}& \lesssim \Vert \nabla v\Vert
_{L^{\infty }}\left( \sum_{j\geq -1}\left( \sum_{|m|\leq 4}2^{(s-1)m}\left(
2^{s(j-m)}\Vert \Delta _{j-m}u\Vert _{K_{p,q}^{\alpha }}\right) \right)
^{r}\right) ^{1/r}  \notag \\
& \lesssim \Vert \nabla v\Vert _{L^{\infty }}\Vert u\Vert
_{BK_{p,q,r}^{\alpha ,s}}.  \label{est:comt_R_4}
\end{align}%
From (\ref{est:comt_R_1})-(\ref{est:comt_R_4}) and (\ref{est:comt_decomp}),
we conclude (\ref{est:comutator_uv_1}). The estimate (\ref%
{est:comutator_uv_2}) follows directly from (\ref{est:comutator_uv_1}) along
with the embedding $BK_{p,q,r}^{\alpha ,s}\hookrightarrow L^{\infty }$ given
in (\ref{l_infty_imersion}).

Now, we turn to estimate (\ref{est:comutator_uv_3}). For that, consider
again the decomposition (\ref{est:comt_decomp}). For $\mathcal{R}_{j}^{1}$,
using the divergence-free condition of $v$, it follows that $\mathcal{R}%
_{j}^{1}=\sum_{j-k\leq 3}\Delta _{j}(\nabla (\Delta _{k}v\tilde{\Delta}%
_{k}u))$. Then, by Bernstein inequality (\ref{est:bernstein_herz_B}), H\"{o}%
lder inequality (\ref{est:holder_herz}) and (\ref{est:l_infty_herz}), we
have that
\begin{equation*}
\Vert \mathcal{R}_{j}^{1}\Vert _{K_{p,q}^{\alpha }}\lesssim \sum_{j-k\leq
3}2^{j}2^{nk/p}\Vert \Delta _{k}v\Vert _{K_{p,q}^{\alpha }}\Vert \tilde{%
\Delta}_{k}u\Vert _{K_{p,q}^{\alpha }}\lesssim \Vert v\Vert _{BK_{p,q,\infty
}^{\alpha ,s}}\sum_{m\leq 3}2^{m}2^{(j-m)(n/p+1-s)}\Vert \Delta _{j-m}u\Vert
_{K_{p,q}^{\alpha }},
\end{equation*}%
where we used the change $m=j-k$. Multiplying both sides by $2^{sj}$, taking
the $\ell ^{r}-$norm and using Minkowski inequality, we arrive at
\begin{equation*}
\left( \sum_{j\geq -1}\left( 2^{(s-1)j}\Vert \mathcal{R}_{j}^{1}\Vert
_{K_{p,q}^{\alpha }}\right)^{r}\right)^{1/r}\hspace{-0.2cm} \lesssim \Vert v\Vert
_{BK_{p,q,r}^{\alpha ,s}}\sum_{m\leq 3}2^{-m\bar{s}_{1}}\left( \sum_{j\geq
-1}\left( 2^{j\bar{s}_{2}}2^{(s-1)(j-m)}\Vert \Delta _{j-m}u\Vert
_{K_{p,q}^{\alpha }}\right) ^{r}\right) ^{1/r} \hspace{-0.1cm},
\end{equation*}%
where $\bar{s}_{1}:=n/p+1-2s$ and $\bar{s}_{2}:=n/p+1-s$. In the case $%
n/p+1-s\leq 0,$ we have $\sup_{j\geq -1}2^{j(n/p+1-s)}<\infty $ and $%
\sum_{m\leq 3}2^{-m(n/p+1-2s)}<\infty $. Thus, it follows that
\begin{equation}
\left( \sum_{j\geq -1}\left( 2^{(s-1)j}\Vert \mathcal{R}_{j}^{1}\Vert
_{K_{p,q}^{\alpha }}\right) ^{r}\right) ^{1/r}\lesssim \Vert v\Vert
_{BK_{p,q,r}^{\alpha ,s}}\Vert u\Vert _{BK_{p,q,r}^{\alpha ,s-1}}.
\label{est:comut_r1_s-1}
\end{equation}

For the parcel $\mathcal{R}_{j}^{2}$, first note that $\mathcal{R}%
_{j}^{2}=\sum_{|j-k|\leq 4}\Delta _{j}(\mathrm{div}\left( S_{k-1}u\Delta
_{k}v\right) )$. Then, proceeding as above, we can estimate
\begin{equation*}
\Vert \mathcal{R}_{j}^{2}\Vert _{K_{p,q}^{\alpha }}\lesssim \sum_{|j-k|\leq
4}2^{j}\Vert S_{k-1}u\Delta _{k}v\Vert _{K_{p,q}^{\alpha }}\lesssim
\sum_{|j-k|\leq 4}2^{j}\Vert S_{k-1}u\Vert _{L^{\infty }}\Vert \Delta
_{k}v\Vert _{K_{p,q}^{\alpha }}.
\end{equation*}%
Using inequality (\ref{est:l_infty_herz}) yields
\begin{equation}
\Vert S_{k-1}u\Vert _{L^{\infty }}\lesssim \sum_{l\leq k-2}2^{nl/p}\Vert
\Delta _{l}u\Vert _{K_{p,q}^{\alpha }}\lesssim \sum_{l\leq
k-2}2^{l(n/p+1-s)}\left( 2^{(s-1)l}\Vert \Delta _{l}u\Vert _{K_{p,q}^{\alpha
}}\right) .  \label{est:sku1}
\end{equation}%
Recall that $\Delta _{l}u=0,$ for all $l\leq -2$. So, for $n/p+1-s\leq 0$,
it follows that
\begin{equation}
\Vert S_{k-1}u\Vert _{L^{\infty }}\lesssim \Vert u\Vert _{BK_{p,q,\infty
}^{\alpha ,s-1}}\lesssim \Vert u\Vert _{BK_{p,q,r}^{\alpha ,s-1}}.
\label{est:sku2}
\end{equation}%
Consequently, making $m=j-k$ \ and using Minkowski inequality lead us to
\begin{align}
\left( \sum_{j\geq -1}\left( 2^{(s-1)j}\Vert \mathcal{R}_{j}^{2}\Vert
_{K_{p,q}^{\alpha }}\right) ^{r}\right) ^{1/r}& \lesssim \Vert u\Vert
_{BK_{p,q,r}^{\alpha ,s-1}}\left( \sum_{j\geq -1}\left(
2^{(s-1)j}2^{j}\sum_{|m|\leq 4}\Vert \Delta _{j-m}v\Vert _{K_{p,q}^{\alpha
}}\right) ^{r}\right) ^{1/r}  \notag \\
& \lesssim \Vert u\Vert _{BK_{p,q,r}^{\alpha ,s-1}}\sum_{|m|\leq
4}2^{sm}\left( \sum_{j\geq -1}\left( 2^{s(j-m)}\Vert \Delta _{j-m}v\Vert
_{K_{p,q}^{\alpha }}\right) ^{r}\right) ^{1/r}  \notag \\
& \lesssim \Vert u\Vert _{BK_{p,q,r}^{\alpha ,s-1}}\Vert v\Vert
_{BK_{p,q,r}^{\alpha ,s}}.  \label{est:comut_r2_s-1}
\end{align}

Next we handle the parcel $\mathcal{R}_{j}^{3}=-\sum_{j-k\leq 2}\Delta
_{k}vS_{k+2}(\Delta _{j}\nabla u)$ as follows
\begin{equation*}
\Vert \mathcal{R}_{j}^{3}\Vert _{K_{p,q}^{\alpha }}\lesssim \sum_{j-k\leq
2}\Vert \Delta _{k}vS_{k+2}(\Delta _{j}\nabla u)\Vert _{K_{p,q}^{\alpha
}}\lesssim \sum_{j-k\leq 2}\Vert \Delta _{k}v\Vert _{K_{p,q}^{\alpha }}\Vert
S_{k+2}(\Delta _{j}\nabla u)\Vert _{L^{\infty }}.
\end{equation*}%
Using the same arguments to get (\ref{est:sku1}) and (\ref{est:sku2}), we
can obtain that
\begin{align*}
\Vert S_{k+2}(\Delta _{j}\nabla u)\Vert _{L^{\infty}}\lesssim 2^{j}\Vert u\Vert _{BK_{p,q,r}^{\alpha ,s-1}}
\end{align*}
for $n/p+1-s\leq 0$%
. Since $s>0,$ proceeding as in (\ref{est:comut_r2_s-1}), it follows that
\begin{equation}
\left( \sum_{j\geq -1}\left( 2^{(s-1)j}\Vert \mathcal{R}_{j}^{3}\Vert
_{K_{p,q}^{\alpha }}\right) ^{r}\right) ^{1/r}\lesssim \Vert u\Vert
_{BK_{p,q,r}^{\alpha ,s-1}}\Vert v\Vert _{BK_{p,q,r}^{\alpha ,s}}.
\label{est:comut_r3_s-1}
\end{equation}

Finally, for $\mathcal{R}_{j}^{4}$, we can estimate
\begin{equation*}
\Vert \mathcal{R}_{j}^{4}\Vert _{K_{p,q}^{\alpha }}\lesssim \sum_{|j-k|\leq
4}2^{-j}\Vert \nabla S_{k-1}v\Vert _{L^{\infty }}\Vert \Delta _{k}\nabla
u\Vert _{K_{p,q}^{\alpha }}.
\end{equation*}
For $n/p+1-s\leq 0$, we also get $\Vert \nabla S_{k-1}v\Vert _{L^{\infty
}}\lesssim \Vert v\Vert _{BK_{p,q,r}^{\alpha ,s}}$. Now, using Bernstein \ref%
{est:bernstein_herz_B} and Minkowski inequality, we obtain
\begin{equation*}
\left( \sum_{j\geq -1}\left( 2^{(s-1)j}\Vert \mathcal{R}_{j}^{4}\Vert
_{K_{p,q}^{\alpha }}\right) ^{r}\right) ^{1/r}\lesssim \Vert v\Vert
_{BK_{p,q,r}^{\alpha ,s}}\left( \sum_{j\geq -1}\left(
2^{(s-1)j}\sum_{|j-k|\leq 4}2^{k-j}\Vert \Delta _{k}u\Vert _{K_{p,q}^{\alpha
}}\right) ^{r}\right) ^{1/r}.
\end{equation*}%
For $s\in \mathbb{R}$, making $m=k-j$, it follows that
\begin{align}
\left( \sum_{j\geq -1}\left( 2^{(s-1)j}\Vert \mathcal{R}_{j}^{4}\Vert
_{K_{p,q}^{\alpha }}\right) ^{r}\right) ^{1/r}& \lesssim \Vert v\Vert
_{BK_{p,q,r}^{\alpha ,s}}\sum_{|m|\leq 4}2^{m(2-s)}\left( \sum_{j\geq
-1}\left( 2^{(s-1)(j+m)}\Vert \Delta _{j+m}u\Vert _{K_{p,q}^{\alpha
}}\right) ^{r}\right) ^{1/r}  \notag \\
& \lesssim \Vert v\Vert _{BK_{p,q,r}^{\alpha ,s}}\Vert u\Vert
_{BK_{p,q,r}^{\alpha ,s-1}}.  \label{est:comut_r4_s-1}
\end{align}%
From (\ref{est:comut_r1_s-1}), (\ref{est:comut_r2_s-1}), (\ref%
{est:comut_r3_s-1}) and (\ref{est:comut_r4_s-1}), we conclude (\ref%
{est:comutator_uv_3}).

\fin

In order to handle the pressure term $\nabla \pi $, we need a suitable
commutator estimate for $\nabla \pi $ in the context of Besov-Herz spaces.
This is the subject of the next lemma.

\begin{lemma}
\label{lem:comutator_pressure} Let $1\leq p<\infty $, $1\leq q,r\leq \infty $%
, $\alpha \geq 0$ and let $a$, $\pi $ stand for the density and pressure,
respectively.

\begin{itemize}
\item[(i)] For $s>0,$ there exists a constant $C>0$ such that
\begin{equation}
\left( \sum_{j\geq -1}\left( 2^{sj}\Vert \lbrack \Delta _{j},a]\nabla \pi
\Vert _{K_{p,q}^{\alpha }}\right) ^{r}\right) ^{1/r}\leq C\left( \Vert
\nabla \pi \Vert _{L^{\infty }}\Vert a\Vert _{BK_{p,q,r}^{\alpha ,s}}+\Vert
a\Vert _{L^{\infty }}\Vert \nabla \pi \Vert _{BK_{p,q,r}^{\alpha ,s}}\right)
\label{est:comt_pi_1}
\end{equation}%
for all $a,\nabla \pi \in BK_{p,q,r}^{\alpha ,s}$ with $a,\nabla \pi \in
L^{\infty }$.

\item[(ii)] Let $\alpha \geq 0$, $s>n/p$ or $s=n/p$ with $r=1$. Then, we
have that
\begin{equation}
\left( \sum_{j\geq -1}\left( 2^{sj}\Vert \lbrack \Delta _{j},a]\nabla \pi
\Vert _{K_{p,q}^{\alpha }}\right) ^{r}\right) ^{1/r}\leq C\Vert a\Vert
_{BK_{p,q,r}^{\alpha ,s}}\Vert \nabla \pi \Vert _{BK_{p,q,r}^{\alpha ,s}},
\label{est:comt_pi_2}
\end{equation}
for all $a,\nabla \pi \in BK_{p,q,r}^{\alpha ,s},$ where $C>0$ is a constant.

\item[(iii)] Let $\alpha \geq 0$, $s>n/p+1$ or $s=n/p+1$ with $r=1$. Then,
we have that
\begin{equation}
\left( \sum_{j\geq -1}\left( 2^{(s-1)j}\Vert \lbrack \Delta _{j},a]\nabla
\pi \Vert _{K_{p,q}^{\alpha }}\right) ^{r}\right) ^{1/r}\leq C\Vert a\Vert
_{BK_{p,q,r}^{\alpha ,s}}\Vert \nabla \pi \Vert _{BK_{p,q,r}^{\alpha ,s-1}},
\label{est:comt_pi_3}
\end{equation}%
for all $a\in BK_{p,q,r}^{\alpha ,s}$ and $\nabla \pi \in BK_{p,q,r}^{\alpha
,s-1},$ where $C>0$ is a constant.
\end{itemize}
\end{lemma}

\textbf{Proof:} In view of the Bony decomposition, we can express
\begin{equation}
\lbrack \Delta _{j},a]\nabla \pi =[\Delta _{j},T_{a}]\nabla \pi +\Delta
_{j}\left( R(a,\nabla \pi )\right) -R(a,\Delta _{j}\nabla \pi )=:\mathcal{A}%
_{j}^{1}+\mathcal{A}_{j}^{2}+\mathcal{A}_{j}^{3}.  \label{est:comt_decomp_pi}
\end{equation}%
First, note that $\Delta _{j}(\Delta _{k}aS_{k+2}\nabla \pi )\equiv 0,$ for $%
|j-k|\geq 8$, and that $\mathcal{A}_{j}^{2}=\sum_{|j-k|\leq 7}\Delta
_{j}\left( \Delta _{k}aS_{k+2}\nabla \pi \right) $. Then, we can use Young
inequality in Lemma \ref{lem:young_herz} and H\"{o}lder inequality (\ref%
{est:holder_herz}) to estimate
\begin{equation*}
\Vert \mathcal{A}_{j}^{2}\Vert _{K_{p,q}^{\alpha }}\lesssim \sum_{|j-k|\leq
7}\Vert \Delta _{k}a\Vert _{K_{p,q}^{\alpha }}\Vert S_{k+2}\nabla \pi \Vert
_{L^{\infty }}\lesssim \Vert \nabla \pi \Vert _{L^{\infty }}\sum_{|j-k|\leq
7}\Vert \Delta _{k}a\Vert _{K_{p,q}^{\alpha }},
\end{equation*}%
which leads us to
\begin{equation}
2^{sj}\Vert \mathcal{A}_{j}^{2}\Vert _{K_{p,q}^{\alpha }}\lesssim \Vert
\nabla \pi \Vert _{L^{\infty }}\sum_{|m|\leq 7}2^{ms}\left( 2^{s(j-m)}\Vert
\Delta _{j-m}a\Vert _{K_{p,q}^{\alpha }}\right) .  \label{est:commut_finite}
\end{equation}%
Applying the $\ell ^{r}$-norm in (\ref{est:commut_finite}) and using
Minkowski inequality, it follows that
\begin{align}
\left( \sum_{j\geq -1}\left( 2^{sj}\Vert \mathcal{A}_{j}^{2}\Vert
_{K_{p,q}^{\alpha }}\right) ^{r}\right) ^{1/r}& \lesssim \Vert \nabla \pi
\Vert _{L^{\infty }}\sum_{|m|\leq 7}2^{ms}\left( \sum_{j\geq -1}\left(
2^{s(j-m)}\Vert \Delta _{j-m}a\Vert _{K_{p,q}^{\alpha }}\right) ^{r}\right)
^{1/r}  \notag \\
& \lesssim \Vert \nabla \pi \Vert _{L^{\infty }}\Vert a\Vert
_{BK_{p,q,r}^{\alpha ,s}}.  \label{est:comt_A_2}
\end{align}%
For $\mathcal{A}_{j}^{3}$, observing that $\mathcal{A}_{j}^{3}=-\sum_{j-k%
\leq 2}\Delta _{k}aS_{k+2}(\Delta _{j}(\nabla \pi ))$, we obtain that
\begin{equation}
\Vert \mathcal{A}_{j}^{3}\Vert _{K_{p,q}^{\alpha }}\lesssim \sum_{j-k\leq
2}\Vert \Delta _{k}a\Vert _{K_{p,q}^{\alpha }}\Vert S_{k+2}\nabla \pi \Vert
_{L^{\infty }}\lesssim \Vert \nabla \pi \Vert _{L^{\infty }}\sum_{j-k\leq
2}\Vert \Delta _{k}a\Vert _{K_{p,q}^{\alpha }}.  \label{A3-aux}
\end{equation}%
Now, recalling $s>0$ and proceeding analogously to the proof of (\ref%
{est:commut_finite}) and (\ref{est:comt_A_2}), estimate (\ref{A3-aux}) leads
us to%
\begin{align}
\left( \sum_{j\geq -1}\left( 2^{sj}\Vert \mathcal{A}_{j}^{3}\Vert
_{K_{p,q}^{\alpha }}\right) ^{r}\right) ^{1/r}& \lesssim \Vert \nabla \pi
\Vert _{L^{\infty }}\left( \sum_{j\geq -1}\left( \sum_{m\leq 2}2^{ms}\left(
2^{s(j-m)}\Vert \Delta _{j-m}a\Vert _{K_{p,q}^{\alpha }}\right) \right)
^{r}\right) ^{1/r}  \notag \\
& \lesssim \Vert \nabla \pi \Vert _{L^{\infty }}\Vert a\Vert
_{BK_{p,q,r}^{\alpha ,s}}.  \label{est:comt_A_3}
\end{align}%
Finally, we estimate the parcel $\mathcal{A}_{j}^{1}$. For this term, we
proceed as follows
\begin{equation*}
\Vert \mathcal{A}_{j}^{1}\Vert _{K_{p,q}^{\alpha }}\lesssim \sum_{|j-k|\leq
4}\Vert \Delta _{j}(S_{k-1}a\Delta _{k}\nabla \pi )\Vert _{K_{p,q}^{\alpha
}}+\Vert S_{k-1}a\Delta _{j}(\Delta _{k}\nabla \pi )\Vert _{K_{p,q}^{\alpha
}}\lesssim \Vert a\Vert _{L^{\infty }}\sum_{|j-k|\leq 4}\Vert \Delta
_{k}\nabla \pi \Vert _{K_{p,q}^{\alpha }},
\end{equation*}%
which implies that%
\begin{equation}
2^{sj}\Vert \mathcal{A}_{j}^{1}\Vert _{K_{p,q}^{\alpha }}\lesssim \Vert
a\Vert _{L^{\infty }}\sum_{|m|\leq 4}2^{sm}\left( 2^{(j-m)s}\Vert \Delta
_{j-m}\nabla \pi \Vert _{K_{p,q}^{\alpha }}\right) .
\label{est:commut:finite_2}
\end{equation}%
Thereby, after applying the $\ell ^{r}$-norm in (\ref{est:commut:finite_2}),
we arrive at
\begin{equation}
\left( \sum_{j\geq -1}\left( 2^{sj}\Vert \mathcal{A}_{j}^{1}\Vert
_{K_{p,q}^{\alpha }}\right) ^{r}\right) ^{1/r}\lesssim \Vert a\Vert
_{L^{\infty }}\Vert \nabla \pi \Vert _{BK_{p,q,r}^{\alpha ,s}}.
\label{est:comt_A_1}
\end{equation}%
Considering (\ref{est:comt_A_2}), (\ref{est:comt_A_3}) and (\ref%
{est:comt_A_1}) in (\ref{est:comt_decomp_pi}), the resulting estimate is (%
\ref{est:comt_pi_1}). Furthermore, by (\ref{est:comt_pi_1}) and the
embedding $BK_{p,q,r}^{\alpha ,s}\hookrightarrow L^{\infty }$ given in (\ref%
{l_infty_imersion}), we get (\ref{est:comt_pi_2}) directly. Also, using $%
BK_{p,q,r}^{\alpha ,s-1}\hookrightarrow L^{\infty }$ and (\ref{est:comt_pi_1}%
) with $s-1,$ estimate (\ref{est:comt_pi_3}) follows naturally.

\fin

\section{Linear estimates} \label{sec:linear_est}

We shall obtain solutions for system (\ref{sist:Euler_1}) as limit of
solutions of approximate linear systems. Thus, in order to develop this
approach, we need to obtain estimates in our setting for suitable linear
problems associated with (\ref{sist:Euler_1}). In this direction, we start
with estimates for the transport equation linked to the density $a.$

\begin{proposition}
\label{prop:transport} Assume that $1\leq p<\infty $, $1\leq q,r\leq \infty $%
, $\alpha \geq 0,$ and $s\geq n/p+1$ with $r=1$ if $s=n/p+1$. Consider $%
a_{0}\in BK_{p,q,r}^{\alpha ,s}$ and a field $u\in BK_{p,q,r}^{\alpha ,s}$
with $\nabla \cdot u=0$ and $\nabla u\in L^{\infty }(\mathbb{R}^{n})$ for $%
T>0$. If $a\in L_{T}^{\infty }(BK_{p,q,r}^{\alpha ,s})$ is a solution of the
transport equation
\begin{equation}
\left\{
\begin{array}{l}
\partial _{t}a+u\cdot \nabla a=0 \\
a(\cdot ,0)=a_{0}%
\end{array}%
,\right. \hspace{0.5cm}(x,t)\in \mathbb{R}^{n}\times \mathbb{R}^{+},
\label{sist:transport}
\end{equation}%
then we have the estimate
\begin{equation}
\Vert a(t)\Vert _{BK_{p,q,r}^{\alpha ,s}}\leq C\left( \Vert a_{0}\Vert
_{BK_{p,q,r}^{\alpha ,s}}+\int_{0}^{t}\Vert a(\tau )\Vert
_{BK_{p,q,r}^{\alpha ,s}}\Vert u(\tau )\Vert _{BK_{p,q,r}^{\alpha
,s}}\;d\tau \right) ,  \label{est:transport}
\end{equation}%
where $C>0$ is a constant. Moreover,
\begin{equation}
\Vert a\Vert _{L_{T}^{\infty }(BK_{p,q,r}^{\alpha ,s})}\leq C\exp \left(
CT\Vert u\Vert _{L_{T}^{\infty }(BK_{p,q,r}^{\alpha ,s})}\right) \Vert
a_{0}\Vert _{BK_{p,q,r}^{\alpha ,s}}.  \label{est:transport_exp}
\end{equation}
\end{proposition}

\textbf{Proof.} Applying the localization $\Delta _{j}$ to (\ref%
{sist:transport}) and using the commutator operator, we obtain that
\begin{equation}
\left\{
\begin{array}{l}
\partial _{t}\Delta _{j}a+u\cdot \nabla \Delta _{j}a=[u\cdot \nabla ,\Delta
_{j}]a \\
\Delta _{j}a(\cdot ,0)=\Delta _{j}a_{0}%
\end{array}%
,\right. \hspace{0.5cm}(x,t)\in \mathbb{R}^{n}\times \mathbb{R}^{+}.
\label{sist:transport_aj}
\end{equation}%
Considering the flow $X$ given in Remark \ref{rem:field_X}, note that
\begin{equation}
\partial _{t}\left( \Delta _{j}a\left( X(y,t),t\right) \right) =\partial
_{t}\Delta _{j}a\left( X_{j}(y,t),t\right) +\left( u\cdot \nabla \Delta
_{j}a\right) \left( X_{j}(y,t),t\right) .  \label{eq:transport_X_aj}
\end{equation}%
Then, it follows from (\ref{sist:transport_aj}) and (\ref{eq:transport_X_aj}%
) that $\partial _{t}\Delta _{j}a\left( X(y,t),t\right) =\left[ u\cdot
\nabla ,\Delta _{j}\right] a\left( X(y,t),t\right) .$ Integrating from $0$
to $t$ yields
\begin{equation}
\Delta _{j}a\left( X(y,t),t\right) =\Delta _{j}a_{0}(y)+\int_{0}^{t}[u\cdot
\nabla ,\Delta _{j}]a\left( X(y,\tau ),\tau \right) \;d\tau ,
\label{aux-int-equa-1}
\end{equation}%
where we have used the initial condition in (\ref{sist:field_X}). This
equation, together with Lemma \ref{lem:field_x}, leads us to
\begin{align}
\Vert \Delta _{j}a(t)\Vert _{K_{p,q}^{\alpha }}& \lesssim \Vert \Delta
_{j}a\left( X(\cdot ,t),t\right) \Vert _{K_{p,q}^{\alpha }}  \notag \\
& \lesssim \Vert \Delta _{j}a_{0}\Vert _{K_{p,q}^{\alpha
}}+\int_{0}^{t}\Vert \lbrack u\cdot \nabla ,\Delta _{j}]a\left( X(\cdot
,\tau ),\tau \right) \Vert _{K_{p,q}^{\alpha }}\;d\tau  \notag \\
& \lesssim \Vert \Delta _{j}a_{0}\Vert _{K_{p,q}^{\alpha
}}+\int_{0}^{t}\Vert \lbrack u\cdot \nabla ,\Delta _{j}]a(\tau )\Vert
_{K_{p,q}^{\alpha }}\;d\tau .  \label{aux-int-equa-2}
\end{align}%
Next, multiplying by $2^{sj}$ and after taking the $\ell ^{r}$-norm (with $%
j\geq -1$) on both sides of estimate (\ref{aux-int-equa-2}), we arrive at
\begin{equation}
\Vert a(t)\Vert _{BK_{p,q,r}^{\alpha ,s}}\leq C\left( \Vert a_{0}\Vert
_{BK_{p,q,r}^{\alpha ,s}}+\int_{0}^{t}\Vert 2^{sj}\Vert \lbrack u\cdot
\nabla ,\Delta _{j}]a(\tau )\Vert _{K_{p,q}^{\alpha }}\Vert _{\ell
^{r}}\;d\tau \right) ,  \label{est:transport_a1}
\end{equation}%
which, along with the commutator estimates in Lemma \ref{lem:comutator_uv},
implies (\ref{est:transport}). Finally, using Gr\"{o}nwall inequality, we
obtain (\ref{est:transport_exp}) directly from (\ref{est:transport}).

\fin

In the sequel we develop estimates for the inhomogeneous linearized Euler
equations in the framework of Besov-Herz spaces.

\begin{proposition}
\label{prop:Euler} Let $1<p<\infty $, $1\leq q,r\leq \infty $, $0\leq \alpha
<n(1-1/p)$ and $s\geq n/p+1$ with $r=1$ if $s=n/p+1$. Consider $u_{0}\in
BK_{p,q,r}^{\alpha ,s}$, a divergence-free vector field $v\in L_{T}^{\infty
}(BK_{p,q,r}^{\alpha ,s})$, $f\in L_{T}^{1}(BK_{p,q,r}^{\alpha ,s})$ and $%
a\in L_{T}^{\infty }(BK_{p,q,r}^{\alpha ,s})$ for $T>0$. Suppose also that
\begin{align*}
(u,\nabla \pi )\in L_{T}^{\infty }(BK_{p,q,r}^{\alpha ,s})\times
L_{T}^{1}(BK_{p,q,r}^{\alpha ,s})
\end{align*}
solves the linearized Euler equations
\begin{equation}
\left\{
\begin{array}{l}
\partial _{t}u+v\cdot \nabla u+(1+a)\nabla \pi =f \\
\mathrm{div}\;u=0 \\
u(\cdot ,0)=u_{0}%
\end{array}%
,\right. \hspace{0.5cm}(x,t)\in \mathbb{R}^{n}\times \mathbb{R}^{+}.
\label{sist:Euler}
\end{equation}%
Then, there holds
\begin{align} \label{est:Euler}
\Vert u(t)\Vert _{BK_{p,q,r}^{\alpha ,s}}& \leq \left( \Vert u_{0}\Vert
_{BK_{p,q,r}^{\alpha ,s}}+\Vert f\Vert _{L_{t}^{1}(BK_{p,q,r}^{\alpha
,s})}+\int_{0}^{t}\left( 1+\Vert a(\tau )\Vert _{BK_{p,q,r}^{\alpha
,s}}\right) \Vert \nabla \pi (\tau )\Vert _{BK_{p,q,r}^{\alpha ,s}}\right.
\notag \\
& \hspace{7.5cm}\left. +\int_{0}^{t}\Vert u(\tau )\Vert _{BK_{p,q,r}^{\alpha
,s}}\Vert v(\tau )\Vert _{BK_{p,q,r}^{\alpha ,s}}\right) ,
\end{align}%
and
\begin{align} \label{est:pressure}
\Vert \nabla \pi \Vert _{L_{t}^{1}(BK_{p,q,r}^{\alpha ,s})}& \leq C\left(
\Vert f\Vert _{L_{t}^{1}(BK_{p,q,r}^{\alpha ,s})}+\int_{0}^{t}\Vert a(\tau
)\Vert _{BK_{p,q,r}^{\alpha ,s}}\Vert \nabla \pi (\tau )\Vert
_{BK_{p,q,r}^{\alpha ,s}}\;d\tau \right.  \notag \\
& \hspace{6.5cm}\left. +\int_{0}^{t}\Vert u(\tau )\Vert _{BK_{p,q,r}^{\alpha
,s}}\Vert v(\tau )\Vert _{BK_{p,q,r}^{\alpha ,s}}\;d\tau \right) .
\end{align}
As a consequence,
\begin{align} \label{est:Euler_exp}
\Vert u\Vert _{L_{T}^{\infty }(BK_{p,q,r}^{\alpha ,s})}+\Vert \nabla \pi
\Vert _{L_{T}^{1}(BK_{p,q,r}^{\alpha ,s})}& \notag\\
&\hspace{-3cm}\leq C\exp \left(
C\int_{0}^{T}\Vert v(\tau )\Vert _{BK_{p,q,r}^{\alpha ,s}}\;d\tau \right)
\\
& \hspace{-1cm}\times \left( \Vert u_{0}\Vert _{BK_{p,q,r}^{\alpha
,s}}+\Vert f\Vert _{L_{T}^{1}(BK_{p,q,r}^{\alpha ,s})}+\Vert a\Vert
_{L_{T}^{\infty }(BK_{p,q,r}^{\alpha ,s})}\Vert \nabla \pi \Vert
_{L_{T}^{1}(BK_{p,q,r}^{\alpha ,s})}\right) .  \notag
\end{align}
\end{proposition}

\textbf{Proof:} Applying $\Delta _{j}$ to system (\ref{sist:Euler}), we
obtain that
\begin{equation*}
\left\{
\begin{array}{l}
\partial _{t}\Delta _{j}u+v\cdot \nabla \Delta _{j}u=[v\cdot \nabla ,\Delta
_{j}]u-\Delta _{j}((1+a)\nabla \pi )+\Delta _{j}f \\
\Delta _{j}u(\cdot ,0)=\Delta _{j}u_{0}%
\end{array}%
,\right. \hspace{0.5cm}(x,t)\in \mathbb{R}^{n}\times \mathbb{R}^{+}.
\end{equation*}%
Considering the flow $X$ given in Remark \ref{rem:field_X} associated with
the field $v$, we can write
\begin{equation*}
\partial _{t}\Delta _{j}u\left( X(y,t),t\right) =\left[ v\cdot \nabla
,\Delta _{j}\right] u\left( X(y,t),t\right) -\Delta _{j}(a\nabla \pi )\left(
X(y,t),t\right) +\Delta _{j}f\left( X(y,t),t\right) .
\end{equation*}%
Integrating from $0$ to $t$, we arrive at
\begin{equation*}
\Delta _{j}u\left( X(y,t),t\right) =\Delta _{j}u_{0}(y)+\int_{0}^{t}\left(
[v\cdot \nabla ,\Delta _{j}]u-\Delta _{j}((1+a)\nabla \pi )+\Delta
_{j}f\right) \left( X(y,\tau ),\tau \right) \;d\tau ,
\end{equation*}%
which, along with Lemma \ref{lem:field_x}, leads us to
\begin{align}
\Vert \Delta _{j}u(t)\Vert _{K_{p,q}^{\alpha }}& \lesssim \Vert \Delta
_{j}u_{0}\Vert _{K_{p,q}^{\alpha }}+\int_{0}^{t}\left( \Vert \lbrack v\cdot
\nabla ,\Delta _{j}]u(\tau )\Vert _{K_{p,q}^{\alpha }}+\Vert \Delta
_{j}\nabla \pi (\tau )\Vert _{K_{p,q}^{\alpha }}\right) d\tau  \notag \\
& \hspace{5.5cm}+\int_{0}^{t}\left( \Vert \Delta _{j}(a\nabla \pi )(\tau
)\Vert _{K_{p,q}^{\alpha }}+\Vert \Delta _{j}f(\tau )\Vert _{K_{p,q}^{\alpha
}}\right) d\tau .  \label{aux-Euler-nonhom-1}
\end{align}%
Next, multiplying both sides of (\ref{aux-Euler-nonhom-1}) by $2^{sj}$ and
afterwards taking the $\ell ^{r}$-norm with $j\geq -1,$ we get the estimate
\begin{align*}
\Vert u(t)\Vert _{BK_{p,q,r}^{\alpha ,s}}& \lesssim \Vert u_{0}\Vert
_{BK_{p,q,r}^{\alpha ,s}}+\Vert \nabla \pi \Vert
_{L_{t}^{1}(BK_{p,q,r}^{\alpha ,s})}+\Vert a\nabla \pi \Vert
_{L_{t}^{1}(BK_{p,q,r}^{\alpha ,s})}+\Vert f\Vert
_{L_{t}^{1}(BK_{p,q,r}^{\alpha ,s})} \\
& \hspace{7cm}+\int_{0}^{t}\Vert 2^{sj}\Vert \lbrack v\cdot \nabla ,\Delta
_{j}]u(\tau )\Vert _{K_{p,q}^{\alpha }}\Vert _{\ell ^{r}}\;d\tau .
\end{align*}%
Since $\Delta _{j}(a\nabla \pi )=a\Delta _{j}\nabla \pi +[\Delta
_{j},a]\nabla \pi $, we have
\begin{align} \label{est:Euler_1}
\Vert u(t)\Vert _{BK_{p,q,r}^{\alpha ,s}}& \lesssim \Vert u_{0}\Vert
_{BK_{p,q,r}^{\alpha ,s}}+\Vert \nabla \pi \Vert
_{L_{t}^{1}(BK_{p,q,r}^{\alpha ,s})}+\Vert f\Vert
_{L_{t}^{1}(BK_{p,q,r}^{\alpha ,s})}+\int_{0}^{t}\Vert 2^{sj}\Vert a\Delta
_{j}(\nabla \pi )(\tau )\Vert _{K_{p,q}^{\alpha }}\Vert _{\ell ^{r}}\;d\tau
 \notag\\
& \hspace{2.5cm}+\int_{0}^{t}\left( \Vert 2^{sj}\Vert \lbrack \Delta
_{j},a]\nabla \pi (\tau )\Vert _{K_{p,q}^{\alpha }}\Vert _{\ell ^{r}}+\Vert
2^{sj}\Vert \lbrack v\cdot \nabla ,\Delta _{j}]u(\tau )\Vert
_{K_{p,q}^{\alpha }}\Vert _{\ell ^{r}}\right) d\tau .
\end{align}%
Now we can conclude (\ref{est:Euler}) by noting that $\Vert 2^{sj}\Vert
a\Delta _{j}(\nabla \pi )(\tau )\Vert _{K_{p,q}^{\alpha }}\Vert _{\ell
^{r}}\lesssim \Vert a(\tau )\Vert _{L^{\infty }}\Vert \nabla \pi (\tau
)\Vert _{BK_{p,q,r}^{\alpha ,s}}$ and using (\ref{est:Euler_1}) together
with Lemmas \ref{lem:comutator_uv} and \ref{lem:comutator_pressure}.

In the sequel we treat with (\ref{est:pressure}). Applying the divergent in (%
\ref{sist:Euler}) and using $\mathrm{div}\;u=0$, we obtain
\begin{equation}
\mathrm{div}(\nabla \pi )=\mathrm{div}\left( f-a\nabla \pi -v\cdot \nabla
u\right) .  \label{equality=div_pressure}
\end{equation}%
Applying $\dot{\Delta}_{j}$ in the equality above, we get $\mathrm{div}\dot{%
\Delta}_{j}(\nabla \pi )=\mathrm{div}\left( \dot{\Delta}_{j}f-\dot{\Delta}%
_{j}(a\nabla \pi )-\dot{\Delta}_{j}(v\cdot \nabla u)\right) $. Before
proceeding, we point out that the Riesz transforms are bounded in $\dot{B}%
K_{p,q,r}^{\alpha ,s}$, as a consequence of their boundedness in Herz spaces
for $-n/p<\alpha <n(1-1/p)$ (see, e.g., \cite{Lucas_Herz, Li_Yang}). Thus,
by Bernstein inequality (\ref{bernstein_bh_homog}), we can estimate
\begin{equation}
\Vert \nabla \pi \Vert _{\dot{B}K_{p,q,r}^{\alpha ,s}}\lesssim \Vert f\Vert
_{BK_{p,q,r}^{\alpha ,s}}+\Vert a\nabla \pi \Vert _{\dot{B}K_{p,q,r}^{\alpha
,s}}+\Vert \mathrm{div}(v\cdot \nabla u)\Vert _{BK_{p,q,r}^{\alpha ,s-1}},
\label{est:pressure_H}
\end{equation}%
since $BK_{p,q,r}^{\alpha ,s}\hookrightarrow \dot{B}K_{p,q,r}^{\alpha ,s}$
for $s>0$. Employing now the product estimates (\ref{est:prod_uv_1}) and (%
\ref{est:prod_uv_3}) given in Lemma \ref{lem:product_estimates}, it follows
that
\begin{align}
\Vert a\nabla \pi \Vert _{\dot{B}K_{p,q,r}^{\alpha ,s}}+\Vert \mathrm{div}%
(v\cdot \nabla u)\Vert _{BK_{p,q,r}^{\alpha ,s-1}}& \notag\\
&\hspace{-2.5cm}\lesssim \Vert a\Vert
_{L^{\infty }}\Vert \nabla \pi \Vert _{\dot{B}K_{p,q,r}^{\alpha ,s}}+\Vert
a\Vert _{\dot{B}K_{p,q,r}^{\alpha ,s}}\Vert \nabla \pi \Vert _{L^{\infty
}}+\Vert u\Vert _{BK_{p,q,r}^{\alpha ,s}}\Vert v\Vert _{BK_{p,q,r}^{\alpha
,s}}  \notag \\
& \hspace{-2.5cm}\lesssim \Vert a\Vert _{BK_{p,q,r}^{\alpha ,s}}\Vert \nabla
\pi \Vert _{\dot{B}K_{p,q,r}^{\alpha ,s}}+\Vert a\Vert _{BK_{p,q,r}^{\alpha
,s}}\Vert \nabla \pi \Vert _{BK_{p,q,r}^{\alpha ,s}}+\Vert u\Vert
_{BK_{p,q,r}^{\alpha ,s}}\Vert v\Vert _{BK_{p,q,r}^{\alpha ,s}},
\label{est:homog_2}
\end{align}%
for $s\geq n/p+1,$ where $r=1$ in the case $s=n/p+1$. Consequently, from (%
\ref{est:pressure_H}) and (\ref{est:homog_2}), we arrive at
\begin{equation} \label{est:pressure_homog}
\Vert \nabla \pi \Vert _{\dot{B}K_{p,q,r}^{\alpha ,s}}\lesssim \Vert f\Vert
_{BK_{p,q,r}^{\alpha ,s}}+\Vert a\Vert _{BK_{p,q,r}^{\alpha ,s}}\left( \Vert
\nabla \pi \Vert _{\dot{B}K_{p,q,r}^{\alpha ,s}}+\Vert \nabla \pi \Vert
_{BK_{p,q,r}^{\alpha ,s}}\right) +\Vert u\Vert _{BK_{p,q,r}^{\alpha
,s}}\Vert v\Vert _{BK_{p,q,r}^{\alpha ,s}}.
\end{equation}%
Moreover, in view of (\ref{equality=div_pressure}), $\dot{H}K_{p,q}^{\alpha
,0}=K_{p,q}^{\alpha }$ and the boundedness of the Riesz transforms in
Sobolev-Herz spaces (similarly to $\dot{B}K_{p,q,r}^{\alpha ,s}$), we can
employ Bernstein inequality (\ref{bernstein_sh_homog}) and the equivalence
of norms in (\ref{def:new_norm_inhomog}) in order to obtain
\begin{align}
\Vert \nabla \pi \Vert _{K_{p,q}^{\alpha }}& =\Vert \nabla \pi \Vert _{\dot{H%
}K_{p,q}^{\alpha ,0}}\lesssim \Vert \mathrm{div}(\nabla \pi )\Vert _{\dot{H}%
K_{p,q}^{\alpha ,-1}}  \notag \\
& \lesssim \Vert \mathrm{div}\;f\Vert _{\dot{H}K_{p,q}^{\alpha ,-1}}+\Vert
\mathrm{div}(a\nabla \pi )\Vert _{\dot{H}K_{p,q}^{\alpha ,-1}}+\Vert \mathrm{%
div}(v\cdot \nabla u)\Vert _{\dot{H}K_{p,q}^{\alpha ,-1}}  \notag \\
& \lesssim \Vert f\Vert _{K_{p,q}^{\alpha }}+\Vert a\nabla \pi \Vert
_{K_{p,q}^{\alpha }}+\Vert v\cdot \nabla u\Vert _{K_{p,q}^{\alpha }}  \notag
\\
& \lesssim \Vert f\Vert _{BK_{p,q,r}^{\alpha ,s}}+\Vert a\Vert
_{BK_{p,q,r}^{\alpha ,s}}\Vert \nabla \pi \Vert _{K_{p,q}^{\alpha }}+\Vert
u\Vert _{BK_{p,q,r}^{\alpha ,s}}\Vert v\Vert _{BK_{p,q,r}^{\alpha ,s}},
\label{est:pressure_K}
\end{align}%
for $s\geq n/p+1,$ where $r=1$ when $s=n/p+1$. Summing (\ref%
{est:pressure_homog}) with (\ref{est:pressure_K}), and recalling (\ref%
{def:new_norm_inhomog}), it follows that
\begin{equation*}
\Vert \nabla \pi \Vert _{BK_{p,q,r}^{\alpha ,s}}\leq C\left( \Vert f\Vert
_{BK_{p,q,r}^{\alpha ,s}}+\Vert a\Vert _{BK_{p,q,r}^{\alpha ,s}}\Vert \nabla
\pi \Vert _{BK_{p,q,r}^{\alpha ,s}}+\Vert u\Vert _{BK_{p,q,r}^{\alpha
,s}}\Vert v\Vert _{BK_{p,q,r}^{\alpha ,s}}\right).
\end{equation*}%
Then, integrating from $0$ to $t$, we conclude (\ref{est:pressure}). In
turn, estimate (\ref{est:Euler_exp}) follows from (\ref{est:Euler}) and (\ref%
{est:pressure}) along with the Gr\"{o}nwall inequality.

\fin

\section{Proof of Theorem \protect\ref{the:Euler}} \label{sec:existence}

This section is divided into six subsections in which we perform basic steps of the proof of Theorem \ref{the:Euler}.

\subsection{Approximate solutions}

We are going to obtain a solution as limit of a sequence constructed via a
interaction process (see, e.g., \cite{Danchin_2010}). For that, let $%
(a^{0},u^{0}):=(a_{0},u_{0})$ and assume that the triple $%
(a^{m},u^{m},\nabla \pi ^{m})$ is well defined for $t>0$ and belong to the
certain suitable spaces that will be given later. Consider $a^{m+1}$ the
solution of the linear transport equation
\begin{equation}
\left\{
\begin{array}{l}
\partial _{t}a^{m+1}+u^{m}\cdot \nabla a^{m+1}=0, \\
a^{m+1}(\cdot ,0)=a_{0}^{m+1}:=S_{m+1}a_{0},%
\end{array}%
\right.  \label{sist:transporte_m}
\end{equation}%
and define $(u^{m+1},\nabla \pi ^{m+1})$ as the solution of the linearized
Euler equation
\begin{equation}
\left\{
\begin{array}{l}
\partial _{t}u^{m+1}+u^{m}\cdot \nabla u^{m+1}+(1+a^{m+1})\nabla \pi
^{m+1}=f, \\
\mathrm{div}\;u^{m+1}=0, \\
u^{m+1}(\cdot ,0)=u_{0}^{m+1}:=S_{m+1}u_{0},%
\end{array}%
\right.  \label{sist:navier_stokes_m}
\end{equation}%
for each $m\in \mathbb{N}_{0}$.

\subsection{Uniform boundedness} \label{section:unif-bound}

In this part we provide uniform bounds for the
approximate solutions $(a^{m},u^{m},\nabla \pi ^{m})$ constructed via the
interaction process (\ref{sist:transporte_m})-(\ref{sist:navier_stokes_m}).

\bigskip Let $X$ $^{m}$ denote the flow associated with $u^{m}$ as in Remark %
\ref{rem:field_X} and let $T>0.$ For every $t\in \lbrack 0,T],$ note that
\begin{equation}
\left\vert (X^{m})^{\pm 1}(y,t)-y)\right\vert \leq \int_{0}^{t}\Vert
u^{m}(\tau )\Vert _{L^{\infty }}\;d\tau \leq CT\Vert u^{m}\Vert
_{L_{T}^{\infty }(BK_{p,q,r}^{\alpha ,s})}\leq \gamma ,  \label{est:Xm}
\end{equation}
provided that $\left\Vert u^{m}\right\Vert _{L_{T}^{\infty
}(BK_{p,q,r}^{\alpha ,s})}\leq (CT)^{-1}\gamma ,$ for all $m\geq 0$, where $%
\gamma >0$ is a constant.

Assuming (\ref{est:Xm}), we can use estimate (\ref{est:transport_exp}) in
Proposition \ref{prop:transport} to obtain
\begin{equation}
\Vert a^{m+1}\Vert _{L_{T}^{\infty }(BK_{p,q,r}^{\alpha ,s})}\leq C\exp
\left( CT\Vert u^{m}\Vert _{L_{T}^{\infty }(BK_{p,q,r}^{\alpha ,s})}\right)
\Vert a_{0}\Vert _{BK_{p,q,r}^{\alpha ,s}}.  \label{est:a.mmais1}
\end{equation}%
Similarly, estimate (\ref{est:Euler_exp}) in Proposition \ref{prop:Euler}
leads us to
\begin{align} \label{est:upmmais1}
\Vert (u^{m+1},\nabla \pi ^{m+1})\Vert _{F_{T}^{s}}& \leq C\exp \left(
CT\Vert u^{m}\Vert _{L_{T}^{\infty }(BK_{p,q,r}^{\alpha ,s})}\right)
 \notag \\
& \hspace{.25cm}\times \left[ \Vert u_{0}\Vert _{BK_{p,q,r}^{\alpha
,s}}  + \Vert f\Vert_{L_{T}^{1}(BK_{p,q,r}^{\alpha ,s})} + \Vert a^{m+1}\Vert _{L_{T}^{\infty }(BK_{p,q,r}^{\alpha ,s})}\Vert
(u^{m+1},\nabla \pi ^{m+1})\Vert _{F_{T}^{s}} \right],
\end{align}
where
\begin{equation}
\Vert (u^{m+1},\nabla \pi ^{m+1})\Vert _{F_{T}^{s}}:=\Vert u^{m+1}\Vert
_{L_{T}^{\infty }(BK_{p,q,r}^{\alpha ,s})}+\Vert \nabla \pi ^{m+1}\Vert
_{L_{T}^{1}(BK_{p,q,r}^{\alpha ,s})}.  \label{norm-FsT}
\end{equation}
Using (\ref{est:a.mmais1}) in (\ref{est:upmmais1}), we arrive at
\begin{align}  \label{est:up.Fst}
\Vert (u^{m+1},\nabla \pi ^{m+1})\Vert _{F_{T}^{s}}&\leq C\exp \left(
CT\Vert u^{m}\Vert _{L_{T}^{\infty }(BK_{p,q,r}^{\alpha ,s})}\right) \left[
\Vert u_{0}\Vert _{BK_{p,q,r}^{\alpha,s}} +\Vert f\Vert
_{L_{T}^{1}(BK_{p,q,r}^{\alpha,s})} \right. \notag \\
&\hspace{1.8cm} + \left. C\exp \left( CT\Vert u^{m}\Vert _{L_{T}^{\infty
}(BK_{p,q,r}^{\alpha,s})}\right) \Vert a_{0}\Vert _{BK_{p,q,r}^{\alpha
,s}}\Vert (u^{m+1},\nabla\pi ^{m+1})\Vert _{F_{T}^{s}}\right] .
\end{align}

Now, in view of the above estimates, we can proceed by induction. For $m=0$,
the following estimate applies
\begin{equation}
\Vert (u^{1},\nabla \pi ^{1})\Vert _{F_{T_{1}}^{s}}\leq C_{0}\left( \Vert
u_{0}\Vert _{BK_{p,q,r}^{\alpha ,s}}+\Vert f\Vert
_{L_{T_{1}}^{1}(BK_{p,q,r}^{\alpha ,s})}\right) ,\text{ for some }C_{0}>0%
\text{ and }T_{1}\in (0,T].  \label{est:induction_1}
\end{equation}%
In fact, by (\ref{est:Xm}) and (\ref{est:up.Fst}), we have that
\begin{align} \label{est:up1}
\Vert (u^{1},\nabla \pi ^{1})\Vert _{F_{T}^{s}}& \leq C\exp \left( CT\Vert
u_{0}\Vert _{BK_{p,q,r}^{\alpha ,s}}\right)\left[ \Vert u_{0}\Vert _{BK_{p,q,r}^{\alpha,s}} + \Vert f\Vert _{L_{T}^{1}(BK_{p,q,r}^{\alpha ,s})} \right.  \notag  \\
& \hspace{4.5cm} + \left. C \Vert a_{0}\Vert _{BK_{p,q,r}^{\alpha ,s}}\exp \left( CT\Vert
u_{0}\Vert _{BK_{p,q,r}^{\alpha ,s}}\right) \Vert (u^{1},\nabla \pi
^{1})\Vert _{F_{T}^{s}} \right] ,
\end{align}%
and
\begin{equation}
\left\vert (X^{0})^{\pm 1}(y,T)-y\right\vert \leq CT\Vert u_{0}\Vert
_{BK_{p,q,r}^{\alpha ,s}}\leq \gamma ,  \label{est:Xm_0}
\end{equation}%
choosing $\gamma $ appropriately. Considering $0<T_{1}\leq T$ such that $%
CT_{1}\Vert u_{0}\Vert _{BK_{p,q,r}^{\alpha ,s}}\leq \mathrm{ln}(2)$, we
have $\exp (CT_{1}\Vert u_{0}\Vert _{BK_{p,q,r}^{\alpha ,s}})\leq 2$
and estimate (\ref{est:Xm_0}) holds with $\gamma =\mathrm{ln}(2)$. In this way, we
can estimate%
\begin{align*}
\Vert (u^{1},\nabla \pi ^{1})\Vert _{F_{T_{1}}^{s}}& \leq 2C\left( \Vert
u_{0}\Vert _{BK_{p,q,r}^{\alpha ,s}} +\Vert f\Vert_{L_{T_{1}}^{1}(BK_{p,q,r}^{\alpha ,s})} + 2C\Vert a_{0}\Vert _{BK_{p,q,r}^{\alpha
,s}}\Vert (u^{1},\nabla \pi ^{1})\Vert _{F_{T_{1}}^{s}}\right) \\
& =2C\left( \Vert u_{0}\Vert _{BK_{p,q,r}^{\alpha ,s}}+\Vert f\Vert
_{L_{T_{1}}^{1}(BK_{p,q,r}^{\alpha ,s})}\right) +4C^{2}\Vert a_{0}\Vert
_{BK_{p,q,r}^{\alpha ,s}}\Vert (u^{1},\nabla \pi ^{1})\Vert _{F_{T_{1}}^{s}}.
\end{align*}%
So, taking $4C^{2}\Vert a_{0}\Vert _{BK_{p,q,r}^{\alpha ,s}}\leq 1/2$, we
get
\begin{equation*}
\Vert (u^{1},\nabla \pi ^{1})\Vert _{F_{T_{1}}^{s}}\leq 4C\left( \Vert
u_{0}\Vert _{BK_{p,q,r}^{\alpha ,s}}+\Vert f\Vert
_{L_{T_{1}}^{1}(BK_{p,q,r}^{\alpha ,s})}\right)
\end{equation*}%
and (\ref{est:induction_1}) follows with $C_{0}=4C$.

Now let $\widetilde{C}>0$, $0<T_{2}\leq T_{1}$ and $\gamma >0$ such that $%
C_{0}\left( \Vert u_{0}\Vert _{BK_{p,q,r}^{\alpha ,s}}+\Vert f\Vert
_{L_{T_{2}}^{1}(BK_{p,q,r}^{\alpha ,s})}\right) \leq \widetilde{C}/2$, $%
CT_{2}\widetilde{C}\leq \gamma $ and $\exp (CT_{2}\widetilde{C})\leq 2$. We
claim that
\begin{equation}
\Vert (u^{m},\nabla \pi ^{m})\Vert _{F_{T_{2}}^{s}}\leq \widetilde{C},%
\hspace{0.5cm}\text{for all }m\geq 1\text{.}  \label{est:boundedness}
\end{equation}%
First, in view of (\ref{est:induction_1}), we have (\ref{est:boundedness})
with $m=1$. Suppose that (\ref{est:boundedness}) is valid for a given $m\in
\mathbb{N}$. Since $CT_{2}\widetilde{C}\leq \gamma ,$ it follows that (\ref%
{est:Xm}) holds for $T_{2}>0$. Then, by (\ref{est:up.Fst}) and the induction
hypothesis, we can estimate
\begin{align}
\Vert (u^{m+1},\nabla \pi ^{m+1})\Vert _{F_{T_{2}}^{s}}& \leq C\exp \left( CT_{2}\widetilde{C}\right) \left[ \Vert
u_{0}\Vert _{BK_{p,q,r}^{\alpha ,s}} + \Vert f\Vert
_{L_{T_{2}}^{1}(BK_{p,q,r}^{\alpha ,s})} \right. \notag \\
& \hspace{4cm} +\left. C\exp \left( CT_{2}\widetilde{C}\right)
\Vert a_{0}\Vert _{BK_{p,q,r}^{\alpha ,s}}\Vert (u^{m+1},\nabla \pi
^{m+1})\Vert _{F_{T_{2}}^{s}} \right]  \notag \\
&\leq 2C\left( \Vert u_{0}\Vert _{BK_{p,q,r}^{\alpha
,s}}+2C\Vert a_{0}\Vert _{BK_{p,q,r}^{\alpha ,s}}\Vert (u^{m+1},\nabla \pi
^{m+1})\Vert _{F_{T_{2}}^{s}}+\Vert f\Vert
_{L_{T_{2}}^{1}(BK_{p,q,r}^{\alpha ,s})}\right) .  \label{est:u_mmais1_1}
\end{align}%
Considering $4C^{2}\Vert a_{0}\Vert _{BK_{p,q,r}^{\alpha ,s}}\leq 1/2$ and
recalling $C_{0}=4C,$ estimate (\ref{est:u_mmais1_1}) implies that
\begin{equation}
\Vert (u^{m+1},\nabla \pi ^{m+1})\Vert _{F_{T_{2}}^{s}}\leq 4C\left( \Vert
u_{0}\Vert _{BK_{p,q,r}^{\alpha ,s}}+\Vert f\Vert
_{L_{T_{2}}^{1}(BK_{p,q,r}^{\alpha ,s})}\right) \leq \widetilde{C}.
\label{est:induction_new_2}
\end{equation}%
Then, it follows from (\ref{est:induction_1}), (\ref{est:boundedness}) and (%
\ref{est:induction_new_2}) that $(u^{m+1},\nabla \pi ^{m+1})$ is uniformly
boundedness with respect to the norm $\Vert \cdot \Vert _{F_{T}^{s}}$
defined in (\ref{norm-FsT}). Furthermore, using the boundedness of $%
\{u^{m}\}_{m\in \mathbb{N}_{0}}$ along with (\ref{est:a.mmais1}) for $%
T_{2}>0 $, we get the uniform boundedness of $\{a^{m}\}_{m\in \mathbb{N}%
_{0}} $ in a similar way.

\begin{remark}
\label{Summary-Boundedness-1} In short, we prove that the approximate
solutions $\{(a^{m},u^{m},\nabla \pi ^{m})\}_{m\in \mathbb{N}_{0}}$ are
uniformly bounded in the space $F_{T_{2}}^{s}$ given by
\begin{equation}
F_{T_{2}}^{s}:=L_{T_{2}}^{\infty }(BK_{p,q,r}^{\alpha ,s})\times
L_{T_{2}}^{\infty }(BK_{p,q,r}^{\alpha ,s})\times
L_{T_{2}}^{1}(BK_{p,q,r}^{\alpha ,s}),  \label{space-aux-1}
\end{equation}%
with the norm
\begin{equation}
\Vert (a^{m},u^{m},\nabla \pi ^{m})\Vert _{F_{T_{2}}^{s}}:=\Vert a^{m}\Vert
_{L_{T_{2}}^{\infty }(BK_{p,q,r}^{\alpha ,s})}+\Vert u^{m}\Vert
_{L_{T_{2}}^{\infty }(BK_{p,q,r}^{\alpha ,s})}+\Vert \nabla \pi ^{m}\Vert
_{L_{T_{2}}^{1}(BK_{p,q,r}^{\alpha ,s})}.  \label{norm-space-aux-1}
\end{equation}
\end{remark}

\subsection{Convergence of the approximation scheme} \label{section:convergence}

In this part we show that the sequence $\{(a^{m},u^{m},\nabla \pi ^{m})\}_{m\in \mathbb{N}%
_{0}}$ is a Cauchy sequence in $F_{T}^{s-1},$ for some $0<T\leq T_{2}$,
where the space $F_{T}^{s-1}$ is defined in (\ref{space-aux-1}). For that,
we consider
\begin{equation*}
\delta a^{m+1}:=a^{m+1}-a^{m};\hspace{0.25cm}\delta u^{m+1}:=u^{m+1}-u^{m};%
\hspace{0.25cm}\delta \pi ^{m+1}:=\pi ^{m+1}-\pi ^{m}.
\end{equation*}%
By (\ref{sist:transporte_m}), the difference $\delta a^{m+1}$ satisfies the
problem
\begin{equation}
\left\{
\begin{array}{l}
\partial _{t}\delta a^{m+1}+u^{m}\cdot \nabla \delta a^{m+1}+\delta
u^{m}\cdot \nabla a^{m}=0, \\
\delta a^{m+1}(\cdot ,0)=\Delta _{m+1}a_{0}.%
\end{array}%
\right.   \label{sist:cauchy_transport}
\end{equation}%
Applying $\Delta _{j}$ on both sides of (\ref{sist:cauchy_transport}), it
follows that
\begin{equation*}
\partial _{t}\Delta _{j}\delta a^{m+1}+u^{m}\cdot \nabla \Delta _{j}\delta
a^{m+1}=[u^{m}\cdot \nabla ,\Delta _{j}]\delta a^{m+1}-\Delta _{j}(\delta
u^{m}\cdot \nabla a^{m}).
\end{equation*}%
Using the same process to prove (\ref{est:transport}) in Proposition \ref%
{prop:transport}, along with estimate (\ref{est:comutator_uv_2}) in Lemma %
\ref{lem:comutator_uv}, (\ref{est:prod_uv_2}) in Lemma \ref%
{lem:product_estimates}, and Bernstein inequality, we arrive at
\begin{align}
\Vert \delta a^{m+1}\Vert _{L_{T}^{\infty }(BK_{p,q,r}^{\alpha ,s-1})}&
\notag \\
& \hspace{-2.25cm}\leq \Vert \Delta _{m+1}a_{0}\Vert _{BK_{p,q,r}^{\alpha
,s-1}}+\int_{0}^{T}\left( \Vert 2^{(s-1)j}\Vert \lbrack u^{m}\cdot \nabla
,\Delta _{j}]\delta a^{m+1}\Vert _{K_{p,q}^{\alpha }}\Vert _{\ell
^{r}}+\Vert \delta u^{m}\cdot \nabla a^{m}\Vert _{BK_{p,q,r}^{\alpha
,s-1}}\right) d\tau   \notag \\
& \hspace{-2.25cm}\leq C\left[ 2^{-(m+1)}\Vert a_{0}\Vert _{BK_{p,q,r}^{\alpha
,s}}\right.   \notag \\
&\hspace{-1cm} + \left. T\left( \Vert \delta a^{m+1}\Vert _{L_{T}^{\infty
}(BK_{p,q,r}^{\alpha ,s})}\Vert u^{m}\Vert _{L_{T}^{\infty
}(BK_{p,q,r}^{\alpha ,s-1})}+\Vert \delta u^{m}\Vert _{L_{T}^{\infty
}(BK_{p,q,r}^{\alpha ,s-1})}\Vert a^{m}\Vert _{L_{T}^{\infty
}(BK_{p,q,r}^{\alpha ,s})}\right) \right] .  \label{est:a.cauchy}
\end{align}

Moreover, using (\ref{sist:navier_stokes_m}), we can see that
\begin{equation}
\left\{
\begin{array}{l}
\partial _{t}\delta u^{m+1}+u^{m}\cdot \nabla \delta u^{m+1}+\delta
u^{m}\cdot \nabla u^{m}+\delta a^{m+1}\nabla \pi ^{m}+a^{m+1}\nabla \delta
\pi ^{m+1}=0, \\
\mathrm{div}\;\delta u^{m+1}=0, \\
\delta u^{m+1}(\cdot ,0)=\Delta _{m+1}u_{0}.%
\end{array}%
\right.   \label{sist:cauchy_Euler}
\end{equation}%
Applying $\Delta _{j}$ on both sides of (\ref{sist:cauchy_Euler}) leads us
to
\begin{equation*}
\partial _{t}\Delta _{j}\delta u^{m+1}+u^{m}\cdot \nabla \Delta _{j}\delta
u^{m+1}=[u^{m}\cdot \nabla ,\Delta _{j}]\delta u^{m+1}-\Delta _{j}\left(
a^{m+1}\nabla \delta \pi ^{m+1}+\delta a^{m+1}\nabla \pi ^{m}+\delta
u^{m}\cdot \nabla u^{m}\right) .
\end{equation*}%
Now, proceeding as in the proof of (\ref{est:Euler}) and (\ref{est:pressure}%
) in Proposition \ref{prop:Euler}, and using (\ref{est:prod_uv_2}), (\ref%
{est:comutator_uv_3}), (\ref{est:comt_pi_3}) and Bernstein inequality, we
obtain that
\begin{align}
\Vert \delta u^{m+1}\Vert _{L_{T}^{\infty }(BK_{p,q,r}^{\alpha ,s-1})}&
+\Vert \nabla \delta \pi ^{m+1}\Vert _{L_{T}^{1}(BK_{p,q,r}^{\alpha ,s-1})}
\notag \\
&\hspace{-0.5cm} \leq C\left[ 2^{-(m+1)}\Vert u_{0}\Vert _{BK_{p,q,r}^{\alpha ,s}}+\Vert
\nabla \delta \pi ^{m+1}\Vert _{L_{T}^{1}(BK_{p,q,r}^{\alpha ,s-1})}\Vert
a^{m+1}\Vert _{L_{T}^{\infty }(BK_{p,q,r}^{\alpha ,s})}\right.   \notag \\
& \hspace{1cm} + T\left( \Vert \delta u^{m+1}\Vert _{L_{T}^{\infty
}(BK_{p,q,r}^{\alpha ,s-1})}+\Vert \delta u^{m}\Vert _{L_{T}^{\infty
}(BK_{p,q,r}^{\alpha ,s-1})}\right) \Vert u^{m}\Vert _{L_{T}^{\infty
}(BK_{p,q,r}^{\alpha ,s})}  \notag \\
& \hspace{5.25cm}+\left. \Vert \delta a^{m+1}\Vert _{L_{T}^{\infty
}(BK_{p,q,r}^{\alpha ,s-1})}\Vert \nabla \pi ^{m}\Vert
_{L_{T}^{1}(BK_{p,q,r}^{\alpha ,s})}\right] .  \label{est:up.cauchy}
\end{align}

So, adding $\Vert \delta a^{m+1}\Vert _{L_{T}^{\infty }(N_{p,q,r}^{s-1})}$
to both sides of (\ref{est:up.cauchy}) and using (\ref{est:a.cauchy}), we
can found $C_{1}>0$ such that
\begin{align*}
\Vert (\delta a^{m+1},\delta u^{m+1},\nabla \delta \pi ^{m+1})\Vert
_{F_{T}^{s-1}}& \\
& \hspace{-4cm}\leq C_{1}\left[ 2^{-(m+1)}\left( \Vert u_{0}\Vert
_{BK_{p,q,r}^{\alpha ,s}}+\Vert a_{0}\Vert _{BK_{p,q,r}^{\alpha ,s}}\right)
+\Vert \nabla \delta \pi ^{m+1}\Vert _{L_{T}^{1}(BK_{p,q,r}^{\alpha
,s-1})}\Vert a^{m+1}\Vert _{L_{T}^{\infty }(BK_{p,q,r}^{\alpha ,s})}\right.
\\
& + T\left( \Vert \delta u^{m+1}\Vert _{L_{T}^{\infty
}(BK_{p,q,r}^{\alpha ,s-1})}+\Vert \delta u^{m}\Vert _{L_{T}^{\infty
}(BK_{p,q,r}^{\alpha ,s-1})}\right) \Vert u^{m}\Vert _{L_{T}^{\infty
}(BK_{p,q,r}^{\alpha ,s})} \\
& \hspace{-3cm}+T\left( \Vert \delta a^{m+1}\Vert _{L_{T}^{\infty
}(BK_{p,q,r}^{\alpha ,s-1})}\Vert u^{m}\Vert _{L_{T}^{\infty
}(BK_{p,q,r}^{\alpha ,s})}+\Vert \delta u^{m}\Vert _{L_{T}^{\infty
}(BK_{p,q,r}^{\alpha ,s-1})}\Vert a^{m}\Vert _{L_{T}^{\infty
}(BK_{p,q,r}^{\alpha ,s})}\right) \\
& \hspace{5.75cm}\left. \times \left( \Vert \nabla \pi ^{m}\Vert
_{L_{T}^{1}(BK_{p,q,r}^{\alpha ,s})}+1\right) \right] .
\end{align*}%
By the uniform boundedness of $\{(a^{m},u^{m},\nabla \pi ^{m})\}_{m\in
\mathbb{N}_{0}}$, we know that
\begin{equation*}
\Vert a^{m}\Vert _{L_{T}^{\infty }(BK_{p,q,r}^{\alpha ,s})}\leq C_{0}\Vert
a_{0}\Vert _{BK_{p,q,r}^{\alpha ,s}}\hspace{0.5cm}\text{and} \hspace{0.5cm}
\Vert u^{m}\Vert _{L_{T}^{\infty }(BK_{p,q,r}^{\alpha ,s})}+\Vert \nabla \pi
^{m}\Vert _{L_{T}^{1}(BK_{p,q,r}^{\alpha ,s})}\leq \widetilde{C},
\end{equation*}%
and then, arranging the terms, we get
\begin{align*}
\Vert (\delta a^{m+1},\delta u^{m+1},\nabla \delta \pi ^{m+1})\Vert
_{F_{T}^{s-1}}& \\
&\hspace{-3cm} \leq C_{1}2^{-(m+1)}\left( \Vert u_{0}\Vert
_{BK_{p,q,r}^{\alpha ,s}}+\Vert a_{0}\Vert _{BK_{p,q,r}^{\alpha ,s}}\right) \\
& \hspace{0.85cm} + 2C_{1}T(\widetilde{C}+1)\left( \widetilde{C}+1+C_{0}\Vert a_{0}\Vert
_{BK_{p,q,r}^{\alpha ,s}}\right) \Vert \delta u^{m}\Vert _{L_{T}^{\infty
}(BK_{p,q,r}^{\alpha ,s-1})} \\
& \hspace{-2.5cm} + C_{1}\left[ C_{0}\Vert a_{0}\Vert _{BK_{p,q,r}^{\alpha
,s}}+2T(\widetilde{C}+1)\left( \widetilde{C}+1+C_{0}\Vert a_{0}\Vert
_{BK_{p,q,r}^{\alpha ,s}}\right) \right]\\
&\hspace{4.6cm} \times \Vert (\delta a^{m+1},\delta
u^{m+1},\nabla \delta \pi ^{m+1})\Vert _{F_{T}^{s-1}}.
\end{align*}

Choosing $0<T_{3}\leq T_{2}$ and considering $\Vert a_{0}\Vert
_{BK_{p,q,r}^{\alpha ,s}}$ such that
\begin{equation*}
C_{1}\left[ C_{0}\Vert a_{0}\Vert _{BK_{p,q,r}^{\alpha ,s}}+4T_{3}(%
\widetilde{C}+1)\left( \widetilde{C}+1+C_{0}\Vert a_{0}\Vert
_{BK_{p,q,r}^{\alpha ,s}}\right) \right] \leq \frac{1}{2},
\end{equation*}%
we have that
\begin{equation*}
\Vert (\delta a^{m+1},\delta u^{m+1},\nabla \delta \pi ^{m+1})\Vert
_{F_{T_{3}}^{s-1}}\leq C_{1}2^{-m}\left( \Vert u_{0}\Vert
_{BK_{p,q,r}^{\alpha ,s}}+\Vert a_{0}\Vert _{BK_{p,q,r}^{\alpha ,s}}\right) +%
\frac{1}{2}\Vert \delta u^{m}\Vert _{L_{T_{3}}^{\infty }(BK_{p,q,r}^{\alpha
,s-1})},
\end{equation*}%
for all $m\geq 1$. Recalling that $u^{0}=u_{0}$ and proceeding recurrently,
we get
\begin{align*}
\Vert (\delta a^{m+1},\delta u^{m+1},\nabla \delta \pi ^{m+1})\Vert_{F_{T_{3}}^{s-1}} & \lesssim \frac{m}{2^{m}}\left( \Vert a_{0}\Vert_{BK_{p,q,r}^{\alpha ,s}}+\Vert u_{0}\Vert _{BK_{p,q,r}^{\alpha ,s}}\right) \\
&\hspace{3.5cm} + \frac{1}{2^{m}}\left( \Vert u^{1}\Vert _{L_{T_{5}}^{\infty}(BK_{p,q,r}^{\alpha ,s-1})}+\Vert u_{0}\Vert _{BK_{p,q,r}^{\alpha,s-1}}\right) \\
&\leq \bar{C}\left( \frac{m+1}{2^{m}}\right),
\end{align*}
which implies the desired Cauchy property of $\{(a^{m},u^{m},\nabla \pi
^{m})\}_{m\in \mathbb{N}_{0}}$ in the Banach space $F_{T_{3}}^{s-1}$.
Further, we denote the corresponding limit by $(a,u,\nabla \pi )$.

\subsection{Existence of solution} \label{section:existence_and_continuity}

According to the subsection \ref{section:convergence}, let $(a,u,\nabla \pi
) $ be the limit of the approximate solutions $\{(a^{m},u^{m},\nabla \pi
^{m})\}_{m\in \mathbb{N}_{0}}$ and $T\in (0,\infty ]$ its existence time.
Using duality and the uniform boundedness of the sequence $\{(a^{m},u^{m},\nabla \pi
^{m})\}_{m\in \mathbb{N}_{0}}$, there exists a subsequence $\{(a^{m_{k}},u^{m_{k}},\nabla \pi ^{m_{k}})\}_{m_{k}\in \mathbb{N}_{0}}$ that converges to the limit $(\widetilde{a},\widetilde{u},\nabla\widetilde{\pi })\in F_{T}^{s}$ in a weak-$\ast $ sense. By the uniqueness
of the weak-$\ast $ limit, it follows that $(a,u,\nabla \pi )=(\widetilde{a},%
\widetilde{u},\nabla \widetilde{\pi })$ and, consequently, $(a,u,\nabla \pi
) $ belongs to the space $F_{T}^{s}$ defined in (\ref{space-aux-1}). Also,
since
\begin{equation*}
\{(a^{m},u^{m})\}_{m\in \mathbb{N}_{0}}\subset C([0,T];BK_{p,q,r}^{\alpha
,s-1})\times C([0,T];BK_{p,q,r}^{\alpha ,s-1})
\end{equation*}%
and $\{(a^{m},u^{m})\}_{m\in \mathbb{N}_{0}}$ converges to $(a,u)$ in $%
L_{T}^{\infty }(BK_{p,q,r}^{\alpha ,s-1})\times
L_{T}^{\infty}(BK_{p,q,r}^{\alpha ,s-1})$, it follows that
\begin{equation}
(a,u)\in C([0,T];BK_{p,q,r}^{\alpha ,s-1})\times C([0,T];BK_{p,q,r}^{\alpha
,s-1}).  \label{eq:limit_2}
\end{equation}%
Using that $(a,u,\nabla \pi )\in F_{T}^{s}$ and (\ref{eq:limit_2}), and
passing the limit in (\ref{sist:transporte_m}) and (\ref%
{sist:navier_stokes_m}), we obtain that $(a,u,\nabla \pi )$ is a solution of
(\ref{sist:Euler_1}). The time-continuity of $a$ and $u$ from $[0,T]$ to $%
BK_{p,q,r}^{\alpha ,s}$ follows by a standard argument by using the fact
that they satisfy (\ref{sist:Euler_1}) in a suitable integral sense combined
with the estimates developed in subsection \ref{section:unif-bound}.

\subsection{Uniqueness of solution} \label{section:uniqueness}

This subsection is devoted to the proof of the uniqueness part in Theorem \ref{the:Euler}. Suppose that
\begin{equation}
(a^{i},u^{i},\nabla \pi ^{i})\in C([0,T_{3}];BK_{p,q,r}^{\alpha ,s})\times
C([0,T_{3}];BK_{p,q,r}^{\alpha ,s})\times L_{T_{3}}^{1}(BK_{p,q,r}^{\alpha
,s}),  \label{unique_space}
\end{equation}%
for $i=1,2$, are solutions of system (\ref{sist:Euler_1}) with the same
initial data $(a_{0},u_{0})$ and external force $f$.

Consider $(\delta a,\delta u,\nabla \delta \pi
):=(a^{2}-a^{1},u^{2}-u^{1},\nabla \pi ^{2}-\nabla \pi ^{1})$ and let $%
(a^{1},u^{1},\nabla \pi ^{1})$ be the solution obtained in subsection \ref%
{section:existence_and_continuity}. We have that $\delta a$ satisfies $%
\delta a(\cdot ,0)=0$ and
\begin{equation*}
\partial _{t}\delta a+u^{2}\cdot \nabla \delta a+\delta u\cdot \nabla
a^{1}=0.
\end{equation*}
Then, proceeding as before (see the proof of (\ref{est:a.cauchy})) and using
$\delta a(\cdot ,0)=0$, for all $0<T\leq T_{3}$, it holds that
\begin{equation}
\Vert \delta a\Vert _{L_{T}^{\infty }(BK_{p,q,r}^{\alpha ,s-1})}\leq
C_{0}T\left( \Vert \delta u\Vert _{L_{T}^{\infty }(BK_{p,q,r}^{\alpha
,s-1})}\Vert a^{1}\Vert _{L_{T}^{\infty }(BK_{p,q,r}^{\alpha ,s})}+\Vert
\delta a\Vert _{L_{T}^{\infty }(BK_{p,q,r}^{\alpha ,s-1})}\Vert u^{2}\Vert
_{L_{T}^{\infty }(BK_{p,q,r}^{\alpha ,s})}\right) .  \label{est:a.unic}
\end{equation}%
Moreover, we have that $(\delta u,\nabla \delta \pi )$ satisfies $\mathrm{div%
}\;\delta u=0$, $\delta u(\cdot ,0)=0$ and
\begin{equation*}
\partial _{t}\delta u+u^{2}\cdot \nabla \delta u=-\delta u\cdot \nabla
u^{1}-\delta a\nabla \pi ^{2}-a^{1}\nabla \delta \pi.
\end{equation*}
Proceeding as in the proof of (\ref{est:up.cauchy}), we obtain that
\begin{align}
\Vert \delta u\Vert _{L_{T}^{\infty }(BK_{p,q,r}^{\alpha ,s-1})}+\Vert
\nabla \delta \pi \Vert _{L_{T}^{1}(BK_{p,q,r}^{\alpha ,s-1})}&  \notag \\
& \hspace{-4cm}\leq C_{1}\left[ \Vert a^{1}\Vert _{L_{T}^{\infty
}(BK_{p,q,r}^{\alpha ,s})}\Vert \nabla \delta \pi \Vert
_{L_{T}^{1}(BK_{p,q,r}^{\alpha ,s-1})}+\Vert \delta a\Vert _{L_{T}^{\infty
}(BK_{p,q,r}^{\alpha ,s-1})}\Vert \nabla \pi ^{2}\Vert
_{L_{T}^{1}(BK_{p,q,r}^{\alpha ,s})}\right.  \notag \\
& \hspace{-0.25cm} +\left. T\Vert \delta u\Vert _{L_{T}^{\infty
}(BK_{p,q,r}^{\alpha ,s-1})}\left( \Vert u^{1}\Vert _{L_{T}^{\infty
}(BK_{p,q,r}^{\alpha ,s})}+\Vert u^{2}\Vert _{L_{T}^{\infty
}(BK_{p,q,r}^{\alpha ,s})}\right) \right] ,  \label{est:up.unic}
\end{align}%
since $\delta u(\cdot ,0)=0$. Putting together (\ref{est:a.unic}) and (\ref%
{est:up.unic}), and recalling (\ref{norm-space-aux-1}), we arrive at {%
\begin{align} \label{est:unic_final_1}
\Vert (\delta a,\delta u,\nabla \delta \pi )\Vert _{F_{T}^{s-1}}&  \notag\\
&\hspace{-1.5cm} \leq
C_{0}C_{1}T\left( \Vert \delta u\Vert _{L_{T}^{\infty }(BK_{p,q,r}^{\alpha
,s-1})}\Vert a^{1}\Vert _{L_{T}^{\infty }(BK_{p,q,r}^{\alpha ,s})}+\Vert
\delta a\Vert _{L_{T}^{\infty }(BK_{p,q,r}^{\alpha ,s-1})}\Vert u^{2}\Vert
_{L_{T}^{\infty }(BK_{p,q,r}^{\alpha ,s})}\right) \\
&\hspace{7.5cm} \times \left( \Vert \nabla \pi ^{2}\Vert _{L_{T}^{1}(BK_{p,q,r}^{\alpha
,s})}+1\right)  \notag \\
&\hspace{-1cm} + C_{1}\Vert a^{1}\Vert _{L_{T}^{\infty }(BK_{p,q,r}^{\alpha ,s})}\Vert
\nabla \delta \pi \Vert _{L_{T}^{1}(BK_{p,q,r}^{\alpha ,s-1})}
\notag \\
&\hspace{2cm} + C_{1}T\Vert \delta u\Vert _{L_{T}^{\infty }(BK_{p,q,r}^{\alpha
,s-1})}\left( \Vert u^{1}\Vert _{L_{T}^{\infty }(BK_{p,q,r}^{\alpha
,s})}+\Vert u^{2}\Vert _{L_{T}^{\infty }(BK_{p,q,r}^{\alpha ,s})}\right) .
\notag
\end{align}
By (\ref{unique_space}), we can take $T>0$ sufficiently small such that
\begin{align}
C_{1}T\left( \Vert u^{1}\Vert _{L_{T}^{\infty }(BK_{p,q,r}^{\alpha,s})}+\Vert u^{2}\Vert _{L_{T}^{\infty }(BK_{p,q,r}^{\alpha ,s})}\right) \leq \frac{1}{4}
\end{align}
and
\begin{equation}
C_{0}C_{1}T\left( \Vert a^{1}\Vert _{L_{T}^{\infty }(BK_{p,q,r}^{\alpha
,s})}+\Vert u^{2}\Vert _{L_{T}^{\infty }(BK_{p,q,r}^{\alpha ,s})}\right)
\left( \Vert \nabla \pi ^{2}\Vert _{L_{T}^{1}(BK_{p,q,r}^{\alpha
,s})}+1\right) \leq \frac{1}{4}.  \label{aux-cond-final-1}
\end{equation}%
Considering these last conditions in (\ref{est:unic_final_1}), we obtain
that
\begin{align*}
\Vert (\delta a,\delta u,\nabla \delta \pi )\Vert _{F_{T}^{s-1}}& \leq \frac{%
1}{4}\left( \Vert \delta u\Vert _{L_{T}^{\infty }(BK_{p,q,r}^{\alpha
,s-1})}+\Vert \delta a\Vert _{L_{T}^{\infty }(BK_{p,q,r}^{\alpha
,s-1})}\right) +\frac{1}{4}\Vert \delta u\Vert _{L_{T}^{\infty
}(BK_{p,q,r}^{\alpha ,s-1})} \\
& \hspace{5.75cm}+C_{1}\Vert a^{1}\Vert _{L_{T}^{\infty }(BK_{p,q,r}^{\alpha
,s})}\Vert \nabla \delta \pi \Vert _{L_{T}^{1}(BK_{p,q,r}^{\alpha ,s-1})}.
\end{align*}%
Since $\Vert a^{1}\Vert _{L_{T}^{\infty }(BK_{p,q,r}^{\alpha ,s})}\leq
C_{3}\Vert a_{0}\Vert _{BK_{p,q,r}^{\alpha ,s}}$, considering $%
C_{1}C_{3}\Vert a_{0}\Vert _{BK_{p,q,r}^{\alpha ,s}}\leq 1/2$ leads us to
\begin{align*}
\Vert (\delta a,\delta u,\nabla \delta \pi )\Vert _{F_{T}^{s-1}}& \leq \frac{%
1}{2}\Vert \delta u\Vert _{L_{T}^{\infty }(BK_{p,q,r}^{\alpha ,s-1})}+\frac{1%
}{4}\Vert \delta a\Vert _{L_{T}^{\infty }(BK_{p,q,r}^{\alpha ,s-1})}+\frac{1%
}{2}\Vert \nabla \delta \pi \Vert _{L_{T}^{1}(BK_{p,q,r}^{\alpha ,s-1})} \\
& \leq \frac{1}{2}\Vert (\delta a,\delta u,\nabla \delta \pi )\Vert
_{F_{T}^{s-1}},
\end{align*}%
which implies the uniqueness property. }

\subsection{Continuous dependence on the data} \label{section:cont-dep}

In this section we treat with the continuity of the
solutions obtained in subsection \ref{section:existence_and_continuity},
with respect to the initial data, in a suitable sense.

Consider $(a_{k},u_{k})$ and $(a,u)$ solutions of system (\ref{sist:Euler_1}%
) with initial data $(a_{0,k},u_{0,k})$ and $(a_{0},u_{0})$, respectively.
Also, assume that $\{(a_{0,k},u_{0,k})\}_{k\in \mathbb{N}}$ is a bounded
sequence in $BK_{p,q,r}^{\alpha ,s}$ such that $a_{0,k}\rightarrow a_{0}$
and $u_{0,k}\rightarrow u_{0}$ in $BK_{p,q,r}^{\alpha ,s-1}$. Since $%
\{(a_{0,k},u_{0,k})\}_{k\in \mathbb{N}}$ is bounded in $BK_{p,q,r}^{\alpha
,s}$, it follows that $\{(a_{k},u_{k})\}_{k\in \mathbb{N}}$ is bounded in $%
L_{T}^{\infty }(BK_{p,q,r}^{\alpha ,s})$ for some $T>0$.

Setting $(\delta a_{k},\delta u_{k})=(a_{k}-a,u_{k}-u)$, we have that
\begin{equation*}
\left\{
\begin{array}{l}
\partial _{t}\delta a_{k}+u_{k}\cdot \nabla \delta a_{k}+\delta u_{k}\cdot
\nabla a=0, \\
\delta a_{k}(\cdot ,0)=\delta a_{0,k}.%
\end{array}%
\right.
\end{equation*}%
For $0<T\leq T_{3}$, proceeding as in the proof of (\ref{est:a.cauchy}), it
holds that
\begin{align*}
\Vert \delta a_{k}\Vert _{L_{T}^{\infty }(BK_{p,q,r}^{\alpha ,s-1})} &\\
&\hspace{-1.5cm} \leq C_{0}\left( \Vert \delta a_{0,k}\Vert _{BK_{p,q,r}^{\alpha
,s-1}}+\int_{0}^{t}\left( \Vert \delta u_{k}\Vert _{BK_{p,q,r}^{\alpha
,s-1}}\Vert a\Vert _{BK_{p,q,r}^{\alpha ,s}}+\Vert \delta a_{k}\Vert
_{BK_{p,q,r}^{\alpha ,s-1}}\Vert u_{k}\Vert _{BK_{p,q,r}^{\alpha ,s}}\right)
d\tau \right) ,
\end{align*}%
and, by Gr\"{o}nwall inequality,
\begin{align}
\Vert \delta a_{k}\Vert _{L_{T}^{\infty }(BK_{p,q,r}^{\alpha ,s-1})}\leq
C_{0}\exp \left( C_{0}T\Vert u_{k}\Vert _{L_{T}^{\infty }(BK_{p,q,r}^{\alpha
,s})}\right)&\notag\\
&\hspace{-1.5cm} \times \left[ \Vert \delta a_{0,k}\Vert _{BK_{p,q,r}^{\alpha
,s-1}}+T\Vert \delta u_{k}\Vert _{L_{T}^{\infty }(BK_{p,q,r}^{\alpha
,s-1})}\Vert a\Vert _{L_{T}^{\infty }(BK_{p,q,r}^{\alpha ,s})}\right]
\hspace{-0.05cm}.  \label{est:dp_a}
\end{align}%
Moreover, we have that $\{(u_{k},\nabla \pi _{k})\}_{k\in \mathbb{N}}$
satisfies
\begin{equation*}
\left\{
\begin{array}{l}
\partial _{t}\delta u_{k}+u_{k}\cdot \nabla \delta u_{k}=-\delta u_{k}\cdot
\nabla u-\delta a_{k}\nabla \pi _{k}-a\nabla \delta \pi _{k}, \\
\mathrm{div}\;\delta u_{k}=0, \\
\delta u(\cdot ,0)=\delta u_{0,k}.%
\end{array}%
\right.
\end{equation*}%
Then, we can proceed as in the proof of (\ref{est:up.cauchy}) to arrive at%
\begin{align}
\Vert \delta u_{k}\Vert _{L_{T}^{\infty }(BK_{p,q,r}^{\alpha ,s-1})}& +\Vert
\nabla \delta \pi _{k}\Vert _{L_{T}^{1}(BK_{p,q,r}^{\alpha ,s-1})}  \notag \\
& \hspace{-2.2cm}\leq C_{1}\left[ \Vert \delta u_{0,k}\Vert
_{BK_{p,q,r}^{\alpha ,s-1}}+\int_{0}^{t}\left( \Vert a\Vert
_{BK_{p,q,r}^{\alpha ,s}}\Vert \nabla \delta \pi _{k}\Vert
_{BK_{p,q,r}^{\alpha ,s-1}}+\Vert \delta a_{k}\Vert _{BK_{p,q,r}^{\alpha
,s-1}}\Vert \nabla \pi _{k}\Vert _{BK_{p,q,r}^{\alpha ,s}}\right) d\tau
\right.  \notag \\
& \hspace{2.2cm}\left. +\;\int_{0}^{t}\Vert \delta u_{k}\Vert
_{L_{T}^{\infty }(BK_{p,q,r}^{\alpha ,s-1})}\left( \Vert u\Vert
_{L_{T}^{\infty }(BK_{p,q,r}^{\alpha ,s})}+\Vert u_{k}\Vert _{L_{T}^{\infty
}(BK_{p,q,r}^{\alpha ,s})}\right) d\tau \right] ,  \notag
\end{align}%
and, by Gr\"{o}nwall inequality,
\begin{align} \label{est:dp_u}
\Vert \delta u_{k}\Vert _{L_{T}^{\infty }(BK_{p,q,r}^{\alpha ,s-1})}& +\Vert
\nabla \delta \pi _{k}\Vert _{L_{T}^{1}(BK_{p,q,r}^{\alpha ,s-1})}  \notag \\
& \hspace{-1cm}\leq C_{1}\exp \left( C_{1}T\left( \Vert u\Vert
_{L_{T}^{\infty }(BK_{p,q,r}^{\alpha ,s})}+\Vert u_{k}\Vert _{L_{T}^{\infty
}(BK_{p,q,r}^{\alpha ,s})}\right) \right)  \notag \\
& \hspace{2cm}\times \left[ \Vert \delta u_{0,k}\Vert
_{BK_{p,q,r}^{\alpha ,s-1}}+\Vert a\Vert _{L_{T}^{\infty
}(BK_{p,q,r}^{\alpha ,s})}\Vert \nabla \delta \pi _{k}\Vert
_{L_{T}^{1}(BK_{p,q,r}^{\alpha ,s-1})} \right. \\
&\hspace{6cm} + \left. \Vert \delta a_{k}\Vert
_{L_{T}^{\infty }(BK_{p,q,r}^{\alpha ,s-1})}\Vert \nabla \pi _{k}\Vert_{L_{T}^{1}(BK_{p,q,r}^{\alpha ,s})}\right]. \notag
\end{align}%
Using now the smallness condition in $\Vert a\Vert _{L_{T}^{\infty
}(BK_{p,q,r}^{\alpha ,s})}$ and the boundedness
\begin{equation}
\Vert u\Vert _{L_{T}^{\infty }(BK_{p,q,r}^{\alpha ,s})}+\Vert u_{k}\Vert
_{L_{T}^{\infty }(BK_{p,q,r}^{\alpha ,s})}+\Vert \nabla \pi _{k}\Vert
_{L_{T}^{1}(BK_{p,q,r}^{\alpha ,s})}\leq \widetilde{C}
\label{est:limitation}
\end{equation}%
in estimate (\ref{est:dp_u}) lead us to
\begin{equation}
\Vert \delta u_{k}\Vert _{L_{T}^{\infty }(BK_{p,q,r}^{\alpha ,s-1})}\leq
2C_{1}\exp (C_{1}\widetilde{C}T)\left[ \Vert \delta u_{0,k}\Vert
_{BK_{p,q,r}^{\alpha ,s-1}}+\widetilde{C}\Vert \delta a_{k}\Vert
_{L_{T}^{\infty }(BK_{p,q,r}^{\alpha ,s-1})}\right] .  \label{est:dc_u}
\end{equation}%
Analogously, from (\ref{est:dp_a}) and (\ref{est:limitation}), we can
estimate
\begin{equation}
\Vert \delta a_{k}\Vert _{L_{T}^{\infty }(BK_{p,q,r}^{\alpha ,s-1})}\leq
C_{0}\exp (C_{0}\widetilde{C}T)\left[ \Vert \delta a_{0,k}\Vert
_{BK_{p,q,r}^{\alpha ,s-1}}+T\Vert \delta u_{k}\Vert _{L_{T}^{\infty
}(BK_{p,q,r}^{\alpha ,s-1})}\Vert a\Vert _{L_{T}^{\infty
}(BK_{p,q,r}^{\alpha ,s})}\right] .  \label{est:dc_a}
\end{equation}%
Thus, it follows from (\ref{est:dc_u}) and (\ref{est:dc_a}) that
\begin{align*}
\Vert \delta a_{k}\Vert _{L_{T}^{\infty }(BK_{p,q,r}^{\alpha ,s-1})}& +\Vert
\delta u_{k}\Vert _{L_{T}^{\infty }(BK_{p,q,r}^{\alpha ,s-1})} \\
& \hspace{-1cm}\leq 2C_{1}\exp (C_{1}\widetilde{C}T)\left[ \Vert \delta
u_{0,k}\Vert _{BK_{p,q,r}^{\alpha ,s-1}}\right. \\
& \hspace{-.3cm} + \left. (\widetilde{C}+1)C_{0}\exp (C_{0}\widetilde{C}T)\left(
\Vert \delta a_{0,k}\Vert _{BK_{p,q,r}^{\alpha ,s-1}}+T\Vert \delta
u_{k}\Vert _{L_{T}^{\infty }(BK_{p,q,r}^{\alpha ,s-1})}\Vert a\Vert
_{L_{T}^{\infty }(BK_{p,q,r}^{\alpha ,s})}\right) \right].
\end{align*}
Again, by either the smallness condition on $\Vert a\Vert _{L_{T}^{\infty
}(BK_{p,q,r}^{\alpha ,s})}$ or $T>0$, we obtain that
\begin{align*}
\Vert \delta a_{k}\Vert _{L_{T}^{\infty }(BK_{p,q,r}^{\alpha ,s-1})}& +\Vert
\delta u_{k}\Vert _{L_{T}^{\infty }(BK_{p,q,r}^{\alpha ,s-1})} \\
& \hspace{-0.55cm}\leq 4C_{1}\exp (C_{1}\widetilde{C}T)\left[ \Vert \delta
u_{0,k}\Vert _{BK_{p,q,r}^{\alpha ,s-1}}+(\widetilde{C}+1)C_{0}\exp (C_{0}%
\widetilde{C}T)\Vert \delta a_{0,k}\Vert _{BK_{p,q,r}^{\alpha ,s-1}}\right]
\\
& \hspace{-0.55cm}\leq 4C_{1}\exp (C_{1}\widetilde{C}T)\left[ 1+(\widetilde{C}%
+1)C_{0}\exp (C_{0}\widetilde{C}T)\right] \left( \Vert \delta a_{0,k}\Vert
_{BK_{p,q,r}^{\alpha ,s-1}}+\Vert \delta u_{0,k}\Vert _{BK_{p,q,r}^{\alpha
,s-1}}\right) .
\end{align*}

Therefore, taking
\begin{equation*}
4C_{1}\exp (C_{1}\widetilde{C}T)\left[ 1+(\widetilde{C}+1)C_{0}\exp (C_{0}%
\widetilde{C}T)\right] \leq \bar{C},
\end{equation*}%
with $\bar{C}>0$ independent of $k$, we get the estimate
\begin{equation}
\Vert \delta a_{k}\Vert _{L_{T}^{\infty }(BK_{p,q,r}^{\alpha ,s-1})}+\Vert
\delta u_{k}\Vert _{L_{T}^{\infty }(BK_{p,q,r}^{\alpha ,s-1})}\leq \bar{C}%
\left( \Vert \delta a_{0,k}\Vert _{BK_{p,q,r}^{\alpha ,s-1}}+\Vert \delta
u_{0,k}\Vert _{BK_{p,q,r}^{\alpha ,s-1}}\right) ,
\label{aux-cont-dep-Lipschitz}
\end{equation}%
from which we conclude the continuous dependence.



\begin{thebibliography}{99}

\bibitem{Beirao-Valli-80-1} {H. Beir\~{a}o da Veiga, A. Valli, On the Euler equations for nonhomogeneous fluids. I. \textit{Rend. Sem. Mat. Univ. Padova} \textbf{63} (1980), 151-168.}

\bibitem{Beirao-Valli-80-2} {H. Beir\~{a}o da Veiga, A. Valli, On the Euler equations for nonhomogeneous fluids. II. \textit{J. Math. Anal. Appl.} \textbf{73} (1980), no. 2, 338-350.}

\bibitem{Beirao-Valli-80-3} {H. Beir\~{a}o da Veiga, A. Valli, Existence of }$C^{\infty }${ \ solutions of the Euler equations for nonhomogeneous fluids. \textit{Comm. Partial Differential Equations} \textbf{5} (1980), no. 2, 95-107.}

\bibitem{Bergh_Lofstrom} {J. Bergh, J. L\"{o}fstr\"{o}m, Interpolation spaces. An introduction. Grundlehren der Mathematischen Wissenschaften \textbf{223}\textit{, Springer-Verlag, Berlin-New York,} 1976.}

\bibitem{Bony} {J.-M. Bony, Calcul symbolique et propagation des singularit\'{e}s pour les \'{e}quations aux d\'{e}riv\'{e}es partielles non lin\'{e}aires. \textit{Ann. Sci. \'{E}cole Norm. Sup. (4)} \textbf{14} (1981), no. 2, 209-246.}

\bibitem{Bourgain_2015} {J. Bourgain, D. Li, Strong ill-posedness of the incompressible Euler equation in borderline Sobolev spaces. \textit{Invent. Math.} \textbf{201} (2015), no. 1, 97-157.}

\bibitem{Chae_2003_2} {D. Chae, J. Lee, Local existence and blow-up criterion of the inhomogeneous Euler equations. \textit{J. Math. Fluid Mech.} \textbf{5} (2003), no. 2, 144-165.}

\bibitem{Chae_2003_3} {D. Chae, On the Euler equations in the critical Triebel-Lizorkin spaces. \textit{Arch. Ration. Mech. Anal.} \textbf{170} (2003), no. 3, 185-210.}

\bibitem{Chae-2004} {D. Chae , Local existence and blow-up criterion for the Euler equations in the Besov spaces, \textit{Asymptotic Analysis} \textbf{38} (2004), 339-358.}

\bibitem{Chemin_2} { J.-Y. Chemin, Fluides parfaits incompressibles. \textit{Ast\'{e}risque}} { \textbf{230} (1995), 177 pp.}

\bibitem{Danchin-2006} { R. Danchin, The inviscid limit for density-dependent incompressible fluids. \textit{Ann. Fac. Sci. Toulouse Math.} \textbf{15} (2006), no. 4, 637-688.}

\bibitem{Danchin_2010} { R. Danchin, On the well-posedness of the incompressible density-dependent Euler equations in the $L^{p}$ framework. \textit{J. Differential Equations} \textbf{248} (2010), no. 8, 2130-2170.}

\bibitem{Danchin-2011} { R. Danchin, F. Fanelli, The well-posedness issue for the density-dependent Euler equations in endpoint Besov spaces. \textit{J. Math. Pures Appl.} \textbf{96} (2011), no. 3, 253-278.}

\bibitem{Hernandez} { E. Hern\'{a}ndez, D. Yang, Interpolation of Herz spaces and applications. \textit{Math. Nachr.} \textbf{205} (1999), 69-87.}

\bibitem{Lucas_BWH} {L.C.F. Ferreira,  J.E. P\'{e}rez-L\'{o}pez, Besov-weak-Herz spaces and global solutions for Navier-Stokes equations. \textit{Pacific Journal of Mathematics} \textbf{296} (2018), no. 1, 57-78}

\bibitem{Lucas_Herz} { L.C.F. Ferreira, J.E. P\'{e}rez-L\'{o}pez, On the theory of Besov-Herz spaces and Euler equations. \textit{Israel J. Math.} \textbf{220} (2017), no. 1, 283-332.}

\bibitem{Cuerva_Herrero} {J. Garc\'{\i}a-Cuerva, M.-J.L. Herrero, A theory of Hardy spaces associated to the Herz spaces. \textit{Proc. London Math. Soc. } \textbf{69} (1994), no. 3, 605-628.}

\bibitem{Grafakos} {L. Grafakos, Modern Fourier analysis. Third edition. Graduate Texts in Mathematics \textbf{250}.\textit{\ Springer, New York,} 2014.}

\bibitem{Grafakos_Li_Yang} { L. Grafakos, X. Li, D. Yang, Bilinear operators on Herz-type Hardy spaces. \textit{Trans. Amer. Math. Soc.} \textbf{350} (1998), no. 3, 1249-1275.}

\bibitem{Herz} { C.S. Herz, Lipschitz spaces and Bernstein's theorem on absolutely convergent Fourier transforms. \textit{J. Math. Mech.} \textbf{18} (1968/69), 283-323.}

\bibitem{Itoh-94} { S. Itoh, Cauchy problem for the Euler equations of a nonhomogeneous ideal incompressible fluid. \textit{J. Korean Math. Soc.} \textbf{31} (1994), no. 3, 367-373.}

\bibitem{Itoh-95} { S. Itoh, Cauchy problem for the Euler equations of a nonhomogeneous ideal incompressible fluid. II. \textit{J. Korean Math. Soc.} \textbf{32} (1995), no. 1, 41-50.}

\bibitem{Itoh-99} { S. Itoh, A. Tani, Solvability of nonstationary problems for nonhomogeneous incompressible fluids and the convergence with vanishing viscosity. \textit{Tokyo J. Math.} \textbf{22} (1999), no. 1, 17-42.}

\bibitem{Johnson} { R. Johnson, Lipschitz spaces, Littlewood-Paley spaces, and convoluteurs. \textit{Proc. London Math. Soc.} \textbf{29} (1974), 127-141.}

\bibitem{Kato} { T. Kato, Strong solutions of the Navier-Stokes equation in Morrey spaces. \textit{Bol. Soc. Brasil. Mat. (N.S.)} \textbf{22} (1992), no. 2, 127-155.}

\bibitem{Lemarie} { P.G. Lemari\'{e}-Rieusset, Recent developments in the Navier-Stokes problem. Chapman \& Hall/CRC Research Notes in Mathematics \textbf{431}, \textit{Chapman \& Hall/CRC, Boca Raton, FL,} 2002. }

\bibitem{Li_Yang} {X. Li, D. Yang, Boundedness of some sublinear operators on Herz spaces. \textit{Illinois J. Math.} \textbf{40} (1996), no. 3, 484-501.}

\bibitem{Takada} { R. Takada, Local existence and blow-up criterion for the Euler equations in Besov spaces of weak type. \textit{J. Evol. Equ.} \textbf{8} (2008), no. 4, 693-725.}

\bibitem{Triebel} { H. Triebel, Theory of function spaces. Monographs in Mathematics \textbf{78}, \textit{Birkh\"{a}user Verlag, Basel,} 1983.}

\bibitem{Tsutsui_2011} { Y. Tsutsui, The Navier-Stokes equations and weak Herz spaces. \textit{Adv. Differential Equations} \textbf{16} (2011), no. 11-12, 1049-1085.}

\bibitem{Valli-Za-88} { A. Valli, W.M. Zaj\k{a}czkowski,  About the motion of nonhomogeneous ideal incompressible fluids. \textit{Nonlinear Anal.} \textbf{12} (1988), no. 1, 43-50.}

\bibitem{Vishik-1998} { M. Vishik, Hydrodynamics in Besov spaces. \textit{Archive for Rational and Mechanical Analysis} \textbf{145} (1998), 197-214.}

\bibitem{Wei_2013} { Z. Wei, Local well-posedness for density-dependent incompressible Euler equations. \textit{Electron. J. Differential Equations} \textbf{146} (2013), 18pp.}

\bibitem{Xu_2005} { J. Xu, Equivalent norms of Herz-type Besov and Triebel-Lizorkin spaces. \textit{J. Funct. Spaces Appl.} \textbf{3} (2005), no. 1, 17-31.}

\bibitem{Xu_Yang_2003} { J. Xu, D. Yang, Applications of Herz-type Triebel-Lizorkin spaces. \textit{Acta Math. Sci. Ser. B (Engl. Ed.)} \textbf{23} (2003), no. 3, 328-338.}

\bibitem{Xu_Yang_2005} { J. Xu, D. Yang, Herz-type Triebel-Lizorkin spaces. I. \textit{Acta Math. Sin. (Engl. Ser.)} \textbf{21} (2005), no. 3, 643-654.}

\bibitem{Zhou_2010} {Y. Zhou, Local well-posedness and regularity criterion for the density dependent incompressible Euler equations. \textit{Nonlinear Anal.} \textbf{73} (2010), no. 3, 750-766.}

\end{thebibliography}
\end{document}